\documentclass[reqno,11pt,a4paper]{amsart}
\usepackage{cases}
\usepackage{amsmath}
\usepackage{graphicx}
 \usepackage{mathrsfs}
\usepackage{bbm}
\usepackage{amsfonts}
\usepackage{amssymb}
\usepackage{color}
\usepackage{multirow}
\usepackage{enumerate}
\usepackage{epstopdf}
\usepackage{cite}
\usepackage{tikz}
\usepackage{algorithmic}
\usetikzlibrary{arrows.meta}
\usepackage{hyperref}
\usepackage{enumitem}
\usepackage{epsfig,subfigure,epstopdf}
\numberwithin{figure}{section}
\usepackage[margin=1.1in]{geometry}
\usepackage{array}
 \usepackage{appendix}
\usepackage{booktabs}
\usepackage{hyperref}
\hypersetup{%
  colorlinks = true,
  citecolor=blue, 
  linkcolor  = black
}

\allowdisplaybreaks

\usepackage{titlesec}
\titleformat{\section}{\vskip10pt\large\bfseries}{\thesection.}{0.5em}{\centering\vspace{5pt}}
\titleformat{\subsection}{\vskip10pt\normalsize\bfseries}{\thesubsection.}{0.5em}{}

\newtheorem{theorem}{Theorem}[section]
\newtheorem{lemma}[theorem]{Lemma}

\newtheorem{remark}[theorem]{Remark}
\newtheorem{example}[theorem]{Example}
\newtheorem{assumption}{Assumption}[section]

\numberwithin{equation}{section}

\def\d{{\rm d}}

\graphicspath{{fig/}}

\title{Dynamic Ritz projection of finite element methods\\ for fluid-structure interaction}
\date{\today}
\author[]{Erik Burman$^*$,\quad Buyang Li$^\dag$,\quad Rong Tang$^\dag$}
\date{}

\thanks{$^*$Department of Mathematics, University College London, Gower street, London WC1E 6BT, United Kingdom.  Email address: e.burman@ucl.ac.uk}

\thanks{$^\dag$Department of Applied Mathematics, The Hong Kong Polytechnic University, Hong
Kong. Email address: buyang.li@polyu.edu.hk and claire.tang@polyu.edu.hk}
\begin{document}
\maketitle

{\small 
\textbf{Abstract.} 
Regardless of the development of various finite element methods for fluid-structure interaction (FSI) problems, optimal-order convergence of finite element discretizations of the FSI problems in the $L^\infty(0,T;L^2)$ norm has not been proved due to the incompatibility between standard Ritz projections and the interface conditions in the FSI problems. To address this issue, we define a dynamic Ritz projection (which satisfies a dynamic interface condition) associated to the FSI problem and study its approximation properties in the $L^\infty(0,T;H^1)$ and $L^\infty(0,T;L^2)$ norms. Existence and uniqueness of the dynamic Ritz projection of the solution, as well as estimates of the error between the solution and its dynamic Ritz projection, are established. By utilizing the established results, we prove optimal-order convergence of finite element methods for the FSI problem in the $L^\infty(0,T;L^2)$ norm. \\

\textbf{Key words.} Fluid-structure interaction, finite element method, dynamic Ritz projection, optimal-order convergence
\\

\textbf{MSC codes.} 65M12, 65M15, 76D05
}

\setlength\abovedisplayskip{4pt}
\setlength\belowdisplayskip{4pt}

\section{Introduction}

Fluid-structure interaction (FSI) problems are pervasive in various fields of science and engineering, where the interaction between fluid flows and structures plays a crucial role. The study of FSI encompasses a wide range of applications, from aerospace engineering, where the behavior of aircraft wings under aerodynamic loads is critical, to biomedical engineering, where the dynamics of blood flow interacting with arterial walls can influence the design of medical devices. In civil engineering, understanding the effects of wind or water flow on buildings and bridges is essential for safe and sustainable infrastructure design. Numerical simulations of the FSI problems is important in providing accurate and efficient solutions to complex problems that would otherwise be intractable, leading to improved designs and enhanced understanding of the physical processes involved. The development of novel numerical methods and rigorous numerical analysis for FSI problems has been an exciting frontier in computational mechanics and computational mathematics.

In this paper, we consider a commonly used FSI model described by the Stokes equations in a fluid region $\Omega_1\subset\mathbb{R}^d$ and the vectorial wave equation in a solid region $\Omega_2\subset\mathbb{R}^d$, coupled on a fixed interface $\Gamma$ which separates the two subdomains, i.e., 
\begin{equation}\label{exact_eqn} 
\left\{ 
\begin{aligned} 
\partial_t u - \nabla \cdot (D(u) - p I ) &= f  &&\mbox{in}\,\,\,  \Omega_1  \\ 
\nabla \cdot u &= 0 &&\mbox{in}\,\,\,  \Omega_1 \\ 
\partial_{tt}\eta - \Delta\eta &= 0 &&\mbox{in}\,\,\, \Omega_2 ,
\end{aligned}
\right.
\end{equation}
which are coupled on the interface via the following continuity conditions:
\begin{equation}\label{exact_interface}
\left\{
\begin{aligned} 
    (D(u) - pI)n &= \partial_n \eta &&{\rm on}~~ \Gamma\\
    u &= \partial_t  \eta &&{\rm on}~~ \Gamma , 
\end{aligned}
\right.
\end{equation}
where $n$ is the unit normal vector on the interface $\Gamma$ pointing to $\Omega_2$. For simplicity, we consider the FSI problem in \eqref{exact_eqn}--\eqref{exact_interface} under the Neumann boundary conditions, i.e., 
\begin{equation}\label{exact_boundary}
\left\{
\begin{aligned} 
 (D(u) - pI)n &= g &&{\rm on}~~ \partial \Omega_1\backslash\Gamma \\
\partial_n\eta &= 0 &&{\rm on}~~ \partial \Omega_2\backslash\Gamma , 
\end{aligned}
\right.
\end{equation}
where $g$ is a given function on $\partial \Omega_1\backslash\Gamma$. 
In addition, the following initial conditions are needed to determine the solution uniquely: 
\begin{equation}\label{exact_initial}
u(0)=u^0,\quad \eta(0)=\eta^0 \quad\mbox{and}\quad \partial_t\eta(0)=w^0 . 
\end{equation}
Then, by utilizing the interface conditions in \eqref{exact_interface}, the solution of \eqref{exact_eqn}--\eqref{exact_interface} can be shown to satisfy the following weak formulation: 
\begin{equation}\label{exact-weak} 
\begin{aligned}
&(\partial_{tt}\eta, \xi )_{\Omega_{2}} + (\nabla \eta, \nabla \xi)_{\Omega_2}   + (\partial_t u, v)_{\Omega_{1}} + (D(u), D(v))_{\Omega_{1}} 
- (p, \nabla \cdot v )_{\Omega_{1}} + (\nabla \cdot u, q)_{\Omega_{1}} \\
&= (f, v)_{\Omega_{1}}   + (g, v)_{\partial \Omega_1\backslash\Gamma} , 
\end{aligned}
\end{equation}
for all $(\xi,v,q)\in H^1(\Omega_2)^d\times H^1(\Omega_1)^d \times L^2(\Omega_1) $ satisfying the interface condition $\xi=v$ on $\Gamma$. 

The weak formulation in \eqref{exact-weak} couples the fluid solution $(u,p)$ and the structure solution $\eta$ together. Correspondingly, numerical methods for solving FSI problems can be classified into monolithic schemes and partitioned schemes, where monolithic schemes refer to those methods which solve a coupled system of fluid and structure equations \cite{Michler-2004,Hubner-2004,Rugonyi-2001,Tezduyar-2007,Hao-Sun-2021,Richter-2017}, and partitioned schemes refer to those methods which are designed to decouple the fluid and structure equations in order to allow using well established packages for the single fluid problem and the single structure problem \cite{Badia-2008,Bukac-2020,Burman-Fernandez-2014,Banks-2014,Seboldt-Bukac-2021}. 

Rigorous numerical analysis of both monolithic and partitioned schemes, as well as velocity-pressure decoupling projection methods, has been studied for FSI problems with a fixed interface in \cite{Ambartsumyan-2018,Astorino-2010,Burman-2022,Fernandez-Mullaert-2016,Fu-Kuang-2022,Tallec-2000}. For all the methods, the errors of numerical solutions using finite elements of polynomial degree $k$ were shown to be 
\begin{align} 
\| u_h^n - u(t_n) \|_{L^2}^2  + \| \eta_h^n - \eta(t_n) \|_{H^1}^2  + \sum_{m=1}^n \tau \| u_h^m - u(t_m) \|_{H^1}^2 
\le C( \tau^{2s} + h^{2k} ) , 
\label{error-0} 
\end{align} 
where $\tau$ and $h$ denote the time stepsize and spatial mesh size in the numerical schemes, and $s$ stands for the order of time discretization in the methods. This error estimate for the velocity approximation is of optimal order in the temporally discrete $L^2(0,T;H^1)$ norm but not optimal with respect to $h$ in standard $L^\infty(0,T;L^2)$ norm. 

Since the early work of Wheeler \cite{Wheeler-1973}, it is well known that the Ritz projection plays an important role in establishing optimal-order $L^\infty(0,T;L^2)$ error estimates of FEMs for parabolic and wave equations. However, the time-dependency of the interface condition $u=\partial_t\eta$ on $\Gamma$ poses significant challenges in formulating a properly defined Ritz projection for the FSI problem with desired optimal-order approximation properties. As a result, the error estimates were all established by comparing the numerical solution with the Lagrange interpolation of the exact solution, or Ritz projection of the single fluid or structure equation. This causes the error estimates in the literature to be of suboptimal order in the $L^\infty(0,T;L^2)$ norm even if the exact solutions are assumed to be sufficiently smooth. 
%

In addition to the standard FSI model in \eqref{exact_eqn}--\eqref{exact_interface}, which describes the interaction between fluid and thick structure, similar error estimates as \eqref{error-0} for FSI with a thin-wall structure were studied in \cite{Bukac-Muha-2016,Fernandez-2013,Oyekole-2018}. Only recently, an optimal-order $L^\infty(0,T;L^2)$ error estimate of FEM for an FSI thin-structure problem was proved in \cite{Li-Sun-Xie-Yu-2024} for a partitioned scheme, by utilizing a properly defined dynamic Ritz projection which satisfies the interface condition. For the FSI thin-structure problem, it is shown that the dual problem of the dynamic Ritz projection is equivalent to a backward parabolic equation of $\xi=(D(\phi) - q I )n$ on the interface $\Gamma$, i.e., 
	\begin{equation}\label{xi-strong-eq-0} 
		- \Delta_{\Gamma} \mathcal{N} {\xi} + \mathcal{N} {\xi} - \partial_t \mathcal{}
		{\xi} = f \,\,\, \text{on}\,\,\, \Gamma\times[0,T), \,\,\,\mbox{with initial condition}\,\,\, {\xi} (T) = 0 , 
	\end{equation}
where $\mathcal{N}: H^{-\frac12}(\Gamma)^d\rightarrow H^{\frac12}(\Gamma)^d$ is the Neumann-to-Dirichlet map associated to the Stokes equations. This special property of FSI thin-structure problem was utilized to prove optimal-order estimates of the error between the exact solution and its dynamic Ritz projection, which are further used to establish optimal-order $L^\infty(0,T;L^2)$ error estimates of FEM for the FSI thin-structure problem with a partitioned scheme. 

However, such results have not been established for any FSI thick-structure problems, where the dual problem associated to the dynamic Ritz projection cannot be converted to a backward parabolic equation on the interface $\Gamma$. Consequently, designing a properly defined dynamic Ritz projection that maintains optimal-order approximation properties for the FSI thick-structure problem remains a significant challenge. This issue will be addressed in the present paper. 

The dynamic Ritz projection for the FSI thick-structure problem, which will be introduced and analyzed in this paper, is considered in a semidiscrete setting (without time discretization). Nevertheless, it can be directly applied to study the convergence of fully discrete FEMs for both monolithic and partitioned schemes. To illustrate, we prove optimal-order convergence of a FEM in the $L^\infty(0,T;L^2)$ norm, with a fundamental monolithic Crank–Nicolson scheme, for the FSI problem in \eqref{exact_eqn}–\eqref{exact_initial}. 

The structure of this paper is organized as follows: In Section \ref{section:main}, we present the assumptions and main results, including the definition and properties of the dynamic Ritz projection for the FSI thick-structure problem, as well as the optimal-order error estimates for the FEM using a monolithic Crank--Nicolson scheme in the $L^\infty(0,T;L^2)$ norm. The proofs of the two main results are presented in Section \ref{section:well-posedness} and Section \ref{section:error}, respectively. Finally, numerical examples are presented in Section 5 to support the theoretical results established in this paper. 


\section{Main results}
\label{section:main}

\subsection{Basic settings: Triangulation and regularity assumptions}
We assume that the domain $\Omega = \Omega_1\cup\Gamma \cup \Omega_2$ is triangulated into shape-regular and quasi-uniform simplices (triangles in 2D or tetrahedra in 3D) such that the triangulations of $\Omega_1$ and $\Omega_2$ share common faces on the interface $\Gamma$. 
For the simplicity of illustration of the ideas without complicating the notations, in order to avoid considering the errors in approximating the geometry of domain $\Omega$ and interface $\Gamma$, we assume that the boundary and interface are fitted by the triangulation exactly without approximation errors. This is true in the following two cases: 
(1) Domains $\Omega_1$ and $\Omega_2$ are both polygons or polyhedra, and the interface $\Gamma$ is flat; (2) isogeometric finite elements are used to fit the geometry of domain $\Omega$. 
Once the construction of the dynamic Ritz projection and proof of optimal-order convergence in the $L^\infty(0,T;L^2)$ norm are realized under this simple setting, our proof can be modified to the more general setting where the triangulations do not fit the boundary and interface exactly, by using isoparametric finite elements and including the errors from approximating domain $\Omega$ by a triangulated domain $\Omega_h$.



For the sake of simplicity in notation, we denote the norms of the Sobolev spaces $H^s(\Omega_j)^d$ and $H^s(\Gamma)^d$ as $\|\cdot \|_{H^s(\Omega_j)}$ and $\|\cdot \|_{H^s(\Gamma)}$, respectively, omitting the dependence on the dimension $d$ in the norms. 

The optimal-order convergence of finite element methods in the $L^\infty(0,T;L^2)$ norm requires the solution of the FSI problem in \eqref{exact_eqn}--\eqref{exact_initial} to be sufficiently smooth. In addition, we need to use the $H^2$ regularity estimates of some dual PDE problems associated to the dynamic Ritz projection introduced in this paper. For this reason, we make the following assumptions on the regularity of solutions to the FSI problem, the Poisson equation and the Stokes equations. 

\begin{assumption}\label{AS0}{\it 
The solution of the FSI problem in \eqref{exact_eqn}--\eqref{exact_initial} is sufficiently smooth. 
}
\end{assumption}

\begin{assumption}[Extension from the boundary]\label{AS-Ext}{\it 
For $(v,r)\in H^{\frac32}(\Gamma)\times H^{1/2}(\partial\Omega_2)$, there exists an extension $w\in H^2(\Omega_2)$ such that $w=v$ on $\Gamma$ and $\partial_nw=r$ on $\partial\Omega_2$, and 
\begin{subequations}\label{H2-w-ext}
\begin{align}
\|w\|_{H^2(\Omega_2)} 
&\le C\|v\|_{H^{\frac32}(\Gamma)} + C\|r\|_{H^{\frac12}(\partial\Omega_2)} .
\end{align}
\end{subequations}
}
\end{assumption}

\begin{assumption}[Regularity of the Poisson equation]\label{AS1}{\it 
The solutions of the Poisson equation 
    \begin{align}\label{Poisson-1}
    \left\{
    \begin{aligned}
    -\Delta w&=f_2 &&\mbox{in}\,\,\,\Omega_2 \\
    \partial_nw&=0 &&\mbox{on}\,\,\,\partial\Omega_2\backslash\Gamma 
    \end{aligned}
    \right.
    \end{align}
have the following regularity properties: 
\begin{subequations}\label{H2-w}
\begin{align}
\|w\|_{H^2(\Omega_2)} 
&\le C\|f_2\|_{L^2(\Omega_2)} + C\|w\|_{H^{\frac32}(\Gamma)} , \\
\|w\|_{H^2(\Omega_2)} 
&\le C\|f_2\|_{L^2(\Omega_2)} + C\|\partial_nw\|_{H^{\frac12}(\Gamma)} .
\end{align}
\end{subequations}
}
\end{assumption}

\begin{assumption}[Regularity of the Stokes equations]\label{AS2}{\it 
The solutions of the Stokes equations  
    \begin{align}\label{Stokes-1}
    \left\{
    \begin{aligned}
    -\nabla \cdot (Dv+qI) + v &= f_1 &&\mbox{in}\,\,\,\Omega_1 \\
    (Dv+qI) n&=0 &&\mbox{on}\,\,\,\partial\Omega_1\backslash\Gamma 
    \end{aligned}
    \right.
    \end{align}
have the following $H^2$ regularity properties: 
\begin{align}\label{H2-v}
\|v\|_{H^2(\Omega_1)} + \|q\|_{H^1(\Omega_1)} 
&\le C\|f_1\|_{L^2(\Omega_1)} + C\|(Dv+qI)n\|_{H^{\frac12}(\Gamma)} .
\end{align}
}
\end{assumption}

For example, the regularity results in Assumptions \ref{AS-Ext}--\ref{AS2} hold in the following scenarios. 
\begin{itemize}
    
        \item \textit{Scenario 1: $\Gamma$ does not intersect $\partial\Omega$.} 
        Assumptions \ref{AS-Ext}--\ref{AS2} hold when the domain $\Omega$ is smooth or convex, the interface $\Gamma$ is smooth, and $\Gamma$ does not intersect $\partial\Omega$. This is the case when a two-dimensional elastic structure is located inside $\Omega$ without touching the boundary $\partial\Omega$, or when a two-dimensional fluid region is enclosed by an elastic structure, or fluid flow between two elastic structures in a periodic channel. In particular, when $\Gamma$ does not intersect $\partial\Omega$, the $H^2$ regularity estimates in \eqref{H2-w} rely on the $H^2$ regularity of the Poisson equation under pure Dirichlet boundary condition and pure Neumann boundary condition. Such results in a smooth or convex domain can be found in \cite[\textsection 2.2.2 and \textsection 3.2.1]{Grisvard}. 
Similarly, \eqref{H2-v} relies the $H^2$ regularity of the Stokes equations under the Neumann boundary condition, and such results in a smooth domain or convex polygon can be found in \cite[Eq. (1.5)]{Dauge-1989}, \cite[Theorem 1.1]{Shibata-Shimizu-2003} and \cite[Figure 9]{Orlt-Sandig-1995}. 
Similarly, when $\Gamma$ does not intersect $\partial\Omega$, the extension result (from boundary to domain) in Assumption \ref{AS-Ext} relies on the extension of pure Dirichlet data and pure Neumann data, and such results can be found in \cite[Theorem 1.5.1.2]{Grisvard}. \medskip

        \item \textit{Scenario 2: $\Gamma$ intersects $\partial\Omega$.}
        If $\Gamma$ intersects $\partial\Omega$, then the mixed boundary value problem of the Poisson equation has corner singularities at the intersection points. In a two-dimensional domain $\Omega$, the extension results in Assumption \ref{AS-Ext} can be changed to (as discussed in \cite[Lemma 3.1]{Li-2022})
\begin{align*}
\|w\|_{H^{2-\epsilon}(\Omega_2)} 
&\le C\|v\|_{H^{\frac32-\epsilon}(\partial\Omega_2)} + C\|r\|_{H^{\frac12-\epsilon}(\partial\Omega_2)} , 
\end{align*}
where $\epsilon$ can be arbitrarily small. 
If $\Gamma$ intersects $\partial\Omega$ with angle $\pi/2$, then the following $H^{2-\epsilon}$ regularity estimates of Poisson and Stokes equations hold (cf. \cite{Dauge-1988}, \cite[Corollary 3.7]{Dauge-1992} and \cite[Figure 9]{Orlt-Sandig-1995}): 
        \begin{subequations}\label{H2-epsilon-w}
\begin{align}
\|w\|_{H^{2-\epsilon}(\Omega_2)} 
&\le C\|f\|_{H^{-\epsilon}(\Omega_2)} + C\|w\|_{H^{\frac32-\epsilon}(\Gamma)} \\
\|w\|_{H^{2-\epsilon}(\Omega_2)} 
&\le C\|f\|_{H^{-\epsilon}(\Omega_2)} + C\|\partial_nw\|_{H^{\frac12-\epsilon}(\Gamma)} ,
\end{align}
\end{subequations}
and 
\begin{align}\label{H2-epsilon-v}
\|v\|_{H^{2-\epsilon}(\Omega_2)} 
&\le C\|g\|_{H^{-\epsilon}(\Omega_2)} + C\|(Dv+qI)n\|_{H^{\frac12-\epsilon}(\Gamma)} ,
\end{align}
where $\epsilon>0$ can be arbitrarily small. 
The error analysis in this paper can be modified, by using the estimates in \eqref{H2-epsilon-w} and \eqref{H2-epsilon-v}, to prove almost optimal-order convergence in the $L^\infty(0,T;L^2)$ norm with an error bound of $O(h^{k+1-\epsilon})$ for sufficiently smooth solutions of the FSI problem, where $\epsilon>0$ can be arbitrarily small. 

\begin{remark}{\upshape 
If $\Gamma$ intersects $\partial\Omega$ with other angles than $\pi/2$, or the elasticity equation (instead of the vectorial wave equation) is used in the structure region, then the regularity of the solution may be below $H^2$ (worse than Scenario 2) and therefore the numerical solution with quasi-uniform mesh may have suboptimal-order convergence. In this case, the order reduction depends on the angle of intersection between $\Gamma$ and $\partial\Omega$ (this angle determines the strongness of corner singularity at the intersection point). A graded mesh towards the corner of the domain needs to be used to improve the convergence to optimal order. This would require analyzing the corner singularities of the dual PDE problems of the dynamic Ritz projection introduced in this paper. This more complex case will be studied in the future by utilizing the dynamic Ritz projection approach developed in this paper. 
}
\end{remark}

\end{itemize}

\subsection{FEM for the FSI problem}
Let $X_h(\Omega_1)\times M_h(\Omega_1)$ be a conforming finite element subspace of $H^1(\Omega_1)^d\times L^2(\Omega_1)$ satisfying the following desired inf-sup condition and approximation properties: 
\begin{align}
&\|q_h\|_{L^2(\Omega_1)} 
\le \sup_{\substack{v_h\in X_h(\Omega_1)\\ v_h\neq 0}} \frac{(q_h,\nabla\cdot v_h)}{\|\nabla v_h\|_{L^2(\Omega_1)} } \quad\mbox{for}\,\,\, q_h\in M_h(\Omega_1) , \\
&\inf_{v_h\in X_h(\Omega_1)}
\big( \|v - v_h\|_{L^2(\Omega_1)} + h\|v - v_h\|_{H^1(\Omega_1)} \big) 
\le C\|v\|_{H^{k+1}(\Omega_1)} h^{k+1} , \\
&\inf_{q_h\in M_h(\Omega_1)}
\|q - q_h\|_{L^2(\Omega_1)} 
\le C\|q\|_{H^{k}(\Omega_1)} h^{k} ,
\end{align}
for $v\in H^{k+1}(\Omega_1)^d$ and $q\in H^{k}(\Omega_1) $. Moreover, we assume that $X_h(\Omega_2)$ is a finite element subspace of $H^1(\Omega_2)^d$
such that 
\begin{align}
&\inf_{w_h\in X_h(\Omega_2)}
\big( \|w - w_h\|_{L^2(\Omega_1)} + h\|w - w_h\|_{H^1(\Omega_2)} \big) 
\le C\|w\|_{H^{k+1}(\Omega_2)} h^{k+1} 
\end{align}
for $w\in H^{k+1}(\Omega_2)^d$. We assume that $ X_h(\Omega_1)$ and  $X_h(\Omega_2)$ share the same nodes on the interface $\Gamma$, and two finite element functions $v_h\in X_h(\Omega_1)$ and $w_h\in X_h(\Omega_2)$ satisfy the condition $v_h=w_h$ on $\Gamma$ if and only if they are equal at the nodes on $\Gamma$. 

The semidiscrete FEM for the fluid-structure interaction problem in \eqref{exact_eqn}--\eqref{exact_boundary} is as follows: Find $(\eta_h, u_h, p_h)\in X_h(\Omega_2) \times X_h(\Omega_1)\times M_h(\Omega_1) $ such that $\partial_t\eta_h = u_h$ on $\Gamma$ and the weak formulation 
\begin{align}\label{numerical_weak_old}
    \begin{aligned}
        &(\partial_{tt}\eta_h, \xi_h )_{\Omega_{2}} + (\nabla \eta_h, \nabla \xi_h)_{\Omega_{2}} + (\partial_t u_h, v_h)_{\Omega_{1}} + (D(u_h), D(v_h))_{\Omega_{1}} \\
        &\quad  - (p_h, \nabla \cdot v_h )_{\Omega_{1}} + (\nabla \cdot u_h, q_h)_{\Omega_{1}} \\
         &= (f, v_h)_{\Omega_{1}} + (g, v)_{\partial \Omega_1\backslash\Gamma} 
    \end{aligned}
\end{align}
holds for all $(\xi_h, v_h, q_h)\in X_h(\Omega_2)\times X_h(\Omega_1)\times M_h(\Omega_1) $ such that $\xi_h=v_h$ on $\Gamma$. In addition, the following initial conditions are used: 
\begin{align}\label{numerical_initial}
\eta_h(0)=R_h^0\eta(0) ,\quad \partial_t\eta_h(0)=I_h\partial_t\eta(0) 
\quad\mbox{and}\quad u_h(0)= I_hu(0) ,
\end{align}
where $R_h^0\eta(0)$ denotes the Ritz projection of $\eta(0)$ with the Dirichlet interface condition $R_h^0\eta(0)=I_h\eta(0)$ on $\Gamma$, determined by the weak formulation
\begin{align}\label{def-Rh0}
(\nabla (\eta(0) - R_h^0\eta(0)), \nabla \xi_h) = 0  
\quad\forall\, \xi_h\in X_h(\Omega_2)\,\,\,\mbox{such that $\xi_h=0$ on $\Gamma$}, 
\end{align}
and $I_h \partial_t\eta(0)$ denotes the Lagrange interpolation of $\partial_t\eta(0)$. 
This specific choice of initial value $\eta_h(0)$ is required in proving optimal-order convergence of the finite element solutions in this paper; see Section \ref{section:error-estimates}.

For the simplicity of illustration, we consider a fully discrete FEM with a Crank--Nicolson method for the time discretization. Let $0=t_0<t_1<\cdots<t_N=T$ be a uniform partition of the time interval $[0,T]$ with uniform stepsize $\tau=T/N$. For any given $(\eta_h^{n}, w_h^{n}, u_h^{n})$, find $(\eta_h^{n+1}, w_h^{n+1}, u_h^{n+1}, p_h^{n+1/2})\in X_h(\Omega_2)\times X_h(\Omega_2)\times X_h(\Omega_1)\times M_h(\Omega_1)$ such that $w_h^{n+1} = u_h^{n+1}$ on $\Gamma$ and the weak formulation 
\begin{subequations}\label{numerical_weak}
\begin{align}
    &(\eta_h^{n+1}-\eta_h^{n})/\tau = (w_h^{n+1}+w_h^{n})/2 \quad \mbox{on}\,\,\,\Omega_2 ,  \\[5pt] 
    &((w_h^{n+1}-w_h^{n})/\tau , \xi_h )_{\Omega_{2}} 
        + (\nabla (\eta_h^{n+1}+\eta_h^{n})/2, \nabla \xi_h)_{\Omega_{2}} 
        + ((u_h^{n+1} - u_h^n)/\tau, v_h)_{\Omega_{1}} \notag\\
        &\quad + (D((u_h^{n+1}+u_h^{n})/2 ) , D(v_h))_{\Omega_{1}}  - (p_h^{n+1/2} , \nabla \cdot v_h )_{\Omega_{1}} 
        + (\nabla \cdot (u_h^{n+1}+u_h^{n})/2, q_h)_{\Omega_{1}} \\
        &  =
         (f(t_{n+1/2}), v_h)_{\Omega_{1}} + (g(t_{n+1/2}), v)_{\partial \Omega_1\backslash\Gamma}  \notag 
\end{align}
\end{subequations}
holds for all $(\xi_h, v_h, q_h)\in X_h(\Omega_2)\times X_h(\Omega_1)\times M_h(\Omega_1)$ such that $\xi_h=v_h$ on $\Gamma$. By using the inf-sup condition of the finite element space, it is easy to show that problem \eqref{numerical_weak} determines a unique numerical solution $(\eta_h^{n+1}, w_h^{n+1}, u_h^{n+1}, p_h^{n+1/2})$ for any given $(\eta_h^{n}, w_h^{n}, u_h^{n})$.

\subsection{Dynamic Ritz projection and error estimates}
In order to prove optimal-order convergence in the $L^2$ norm for the error between the numerical solution and the exact solution, we shall define a ``dynamic Ritz projection'' (with a dynamic interface condition) for the fluid-structure interaction problem in order to fit the second interface condition in \eqref{exact_interface}. The idea is to find $(R_h \eta, R_h u, R_h p) \in X_h(\Omega_{2}) \times X_h(\Omega_{1})\times M_h(\Omega_{1})$ 
satisfying the weak formulation 
\begin{subequations}\label{weak_Ritz}
    \begin{align}
        (\nabla R_h \eta, \nabla \xi_h)_{\Omega_{2}} &- (\nabla \eta, \nabla \xi_h)_{\Omega_2} + (D(R_h u), D(v_h))_{\Omega_{1}} - (D(u), D(v_h))_{\Omega_1} \notag\\
         &+ (R_h u, v_h)_{\Omega_{1}} - (u, v_h)_{\Omega_1 } - (R_h p, \nabla \cdot v_h )_{\Omega_{1}} + (p, \nabla \cdot v_h)_{\Omega_1} \notag\\
         &
         + (\nabla \cdot R_h u, q_h)_{\Omega_{1}} - (\nabla \cdot u, q_h)_{\Omega_1}= 0 
         \label{weak_Ritz_Omega}\\
         \partial_t R_h \eta  &= R_h u \quad\mbox{on}\,\,\,\Gamma
         \label{weak_Ritz_Gamma}\\
         R_h \eta (0) &= I_h \eta(0) \quad\mbox{on}\,\,\,\Gamma
         \label{weak_Ritz_t=0}
    \end{align}
\end{subequations} 
for all $(\xi_h, v_h, q_h) \in X_h(\Omega_{2}) \times X_h(\Omega_{1})\times M_h(\Omega_{1})$ such that $\xi_h = v_h$ on interface $\Gamma$, 
where $I_h$ denotes the Lagrange interpolation operator and the initial condition in \eqref{weak_Ritz_t=0} is imposed for the ordinary differential equation (ODE) in \eqref{weak_Ritz_Gamma}. 

Existence and uniqueness of solutions to \eqref{weak_Ritz}, and approximation properties of the dynamic Ritz projection defined in \eqref{weak_Ritz} are proved in this paper. The results are summarized in the following theorem.

\begin{theorem}[Dynamic Ritz projection]\label{THM-Ritz}{\it 
Under Assumption \ref{AS0}--\ref{AS2}, problem \eqref{weak_Ritz} determines a unique solution $(R_h \eta, R_h u, R_h p) \in X_h(\Omega_{2}) \times X_h(\Omega_{1})\times M_h(\Omega_{1})$, called the dynamic Ritz projection of the solution. Moreover, this dynamic Ritz projection has the following approximation properties:
\begin{align} 
\sum_{j=0}^3 
\big( \|\partial_t^j ( R_h u -  u)\|_{L^\infty(0,T;H^1(\Omega_1))}
+ \| \partial_t^j ( R_h \eta -  \eta)\|_{L^\infty(0,T;H^1(\Omega_2))} \big) 
&\le Ch^k , \label{Ritz-approx1}\\
\sum_{j=0}^2
\big( \|\partial_t^j ( R_h u -  u)\|_{L^\infty(0,T;L^2(\Omega_1))}
+ \| \partial_t^j ( R_h \eta -  \eta)\|_{L^\infty(0,T;L^2(\Omega_2))} \big) 
&\le Ch^{k+1} . \label{Ritz-approx2}
\end{align}
}
\end{theorem}

By utilizing the results of Theorem \ref{THM-Ritz}, we prove the following result on the optimal-order convergence of finite element solutions to the FSI problem. 

\begin{theorem}[Optimal-order convergence]\label{THM-Error}{\it 
Under Assumption \ref{AS0}--\ref{AS2}, the numerical solution determined by \eqref{numerical_weak} has the following error bound: 
\begin{align}
   \sup_{1\le n\le N} \| u_h^n - u(t_n)\|_{L^2(\Omega_1)} 
   + \sup_{1\le n\le N} \|  \eta_h^n- \eta(t_n)\|_{H^1(\Omega_2)}   \le C(\tau^2 +  h^{k+1}).
\end{align}
}
\end{theorem}

The rest of this paper is devoted to the proofs of Theorem \ref{THM-Ritz} and Theorem \ref{THM-Error}.

\section{Proof of Theorem \ref{THM-Ritz}}\label{section:well-posedness}

We prove existence and uniqueness of the dynamic Ritz projection in Section \ref{section:ODE}, and prove the approximation properties \eqref{Ritz-approx1}--\eqref{Ritz-approx2} in the remaining subsections. 

\subsection{Existence and uniqueness of the dynamic Ritz projection}\label{section:ODE}
Existence and uniqueness of solutions to \eqref{weak_Ritz} can be shown as follows, by considering it as an ODE on interface $\Gamma$. We only need to show that, at any given time $t\in[0,T]$, $R_hu|_\Gamma$ is uniquely determined by $R_h \eta|_\Gamma$ (with linear dependence on $R_h\eta|_\Gamma$) through the weak formulation in \eqref{weak_Ritz_Omega}. Then \eqref{weak_Ritz} can be reformulated as an ODE in the form of $\partial_tR_h\eta = JR_h\eta$, where $J$ is some linear operator on $X_h(\Gamma) = \{w_h|_{\Gamma}: w_h\in X_h(\Omega_{2}) \}$. Such a linear ODE must have a unique solution. 

In fact, by choosing $v_h = 0 $ and $q_h=0$ in $\Bar{\Omega}_{1}$ and $\xi_h = 0 $ on $\Gamma$, we can obtain the following relation from \eqref{weak_Ritz_Omega}: 
\begin{equation}\label{eta_indep}
    (\nabla (R_h \eta - \eta), \nabla \xi_h)_{\Omega_{2}} = 0 \quad\forall\, \xi_h\in X_h(\Omega_{2})\,\,\,\mbox{such that $\xi_h=0$ on $\Gamma$}. 
\end{equation}
Therefore, if $R_h \eta|_\Gamma$ is given then \eqref{eta_indep} determines the value of $R_h \eta$ in $\Omega_2$ as the finite element solution of the Poisson equation under the mixed boundary condition (i.e., Neumann boundary condition on $\partial\Omega_2\backslash\Gamma$ and Dirichlet boundary condition on $\Gamma$). 

After the value of $R_h \eta$ in $\Omega_2$ is determined, the value of $R_hu$ in $\Omega_1$ can be determined by \eqref{weak_Ritz_Omega} as follows. For any admissible test functions $(v_h, q_h)$, we denote by $\xi_h = \text{Ext}(v_h|_{\Gamma})$ an extension of $v_h|_{\Gamma}$ to the domain $\Omega_{2}$ and regard the term $(\nabla R_h \eta, \nabla \xi_h)_{\Omega_{2}} - (\nabla \eta, \nabla \xi_h)_{\Omega_2}$ in \eqref{weak_Ritz_Omega} as an inhomogeneous right-hand side. Consequently, we solve the Stokes system for $(R_h u, R_h p)$ based on the solution $R_h \eta$ in $\Omega_{2}$. This enables us to determine the interface value of $R_h u$, i.e., $R_h u|_{\Gamma}$. Clearly, the map from $R_h \eta|_{\Gamma}$ to $R_h u|_{\Gamma}$ is linear. 

This proves the existence and uniqueness of the solution to \eqref{weak_Ritz}.\hfill\qed

\subsection{Estimates of $\|\nabla (R_h \eta - \eta)\|_{L^\infty(0,T;L^2(\Omega_2))}$ and $\|R_h u - u\|_{L^2(0,T;H^1(\Omega_1))}$}

By choosing test functions $\xi_h = \partial_t (R_h \eta - I_h \eta)$, $v_h = R_h u - I_h u$ and $q_h = R_h p - I_h p$ in the weak formulation \eqref{weak_Ritz}, which satisfy the interface condition $\xi_h = v_h$ on $\Gamma$, we have 
\begin{equation}\label{L2H1-eqn}
    \begin{aligned}
        & \frac{1}{2}\frac{\d}{\d t}\|\nabla(R_h \eta - I_h \eta)\|_{L^2(\Omega_{2})}^2 + \|D(R_h u - I_h u)\|_{L^2(\Omega_{1})}^2 + \|R_h u - I_h u\|_{L^2(\Omega_{1})}^2  \\
         =&\, (\nabla(\eta - I_h \eta), \nabla \xi_h)_{\Omega_{2}} + (D(u - I_h u), D v_h)_{\Omega_{1}} + (u - I_h u, v_h)_{\Omega_{1}}\\
        & + (I_h p - p, \nabla \cdot v_h)_{\Omega_{1}} + (\nabla \cdot(u - I_h u), q_h)_{\Omega_{1}}\\
        =&\!: A_1 + A_2 + A_3 + A_4 + A_5.
    \end{aligned}
\end{equation}
Integrating \eqref{L2H1-eqn} from $0$ to $t$ with respect to time, and using Korn's inequality (see \cite[Theorem 1.1-2]{Ciarlet-1997}) 
$$
\|v\|_{H^1(\Omega_1)}
\le C(\|v\|_{L^2(\Omega_1)} + \|D(v)\|_{L^2(\Omega_1)}) ,
$$
we obtain
\begin{align}\label{L2H1-left}
  &\frac{1}{2}\|\nabla(R_h \eta - I_h \eta)(t)\|_{L^2(\Omega_{2})}^2 - \frac{1}{2}\|\nabla(R_h \eta - I_h \eta)(0)\|_{L^2(\Omega_{2})}^2 + c_0 \|R_h u - u\|_{L^2(0,t;H^1(\Omega_1))}^2 \notag\\
  &\le \int_0^t [ A_1(s) + A_2(s)  + A_3(s)  + A_4(s)  + A_5(s) ]\d s . 
\end{align}
Note that 
\begin{align}
  &\big|\int_0 ^ t A_1(s) \d s\big| \notag \\
  & = \Big|\int_0^t \partial_t (\nabla (\eta- I_h \eta), \nabla (R_h \eta - I_h \eta))_{\Omega_2} \d s- \int_0 ^t (\nabla (\partial_t \eta - I_h \partial_t \eta), \nabla (R_h \eta - I_h \eta))_{\Omega_2} \d s\Big| \notag\\
        & \leq \Big|(\nabla (\eta- I_h \eta)(t), \nabla (R_h \eta - I_h \eta)(t))_{\Omega_2} - (\nabla (\eta- I_h \eta)(0), \nabla (R_h \eta - I_h \eta)(0))_{\Omega_2}\Big| \notag\\
        & \quad \, + \Big|\int_0 ^t (\nabla (\partial_t \eta - I_h \partial_t \eta), \nabla (R_h \eta - I_h \eta))_{\Omega_2} \d s \Big|\\
        & \leq Ch^k \|\nabla(R_h \eta - I_h \eta)(t)\|_{L^2(\Omega_{2})} + Ch^k \|\nabla(R_h \eta - I_h \eta)(0)\|_{L^2(\Omega_{2})} \notag\\
        & \quad \, + Ch^k \|\nabla(R_h \eta - I_h \eta)\|_{L^\infty(0,t;L^2(\Omega_2))} \notag\\
        & \leq \epsilon \|\nabla(R_h \eta - I_h \eta)(t)\|_{L^2(\Omega_2)}^2 + \epsilon \|\nabla(R_h \eta - I_h \eta)(0)\|_{L^2(\Omega_2)}^2 \notag\\
        &\quad\, + \epsilon \|\nabla(R_h \eta - I_h \eta)\|_{L^\infty(0,t;L^2(\Omega_2))}^2 +  C_\epsilon h^{2k}
\end{align}
and 
\begin{equation}
    \begin{aligned}
        &\Big | \int_0^t [  A_2(s)  + A_3(s)  + A_4(s)  + A_5(s) ]\d s \Big| \\
        &\leq Ch^k \|R_h u - u\|_{L^2(0,t;H^1(\Omega_1))} + Ch^k \|R_h p - I_h p\|_{L^2(0,t;L^2(\Omega_1))} \\
        &\leq \epsilon \|R_h u - u\|_{L^2(0,t;H^1(\Omega_1))}^2 + \epsilon \|R_h p - I_h p\|_{L^2(0,t;L^2(\Omega_1))}^2 + C_\epsilon h^{2k}.
    \end{aligned}
\end{equation}
The subsequent step involves estimating $\|R_h p - I_h p(t)\|_{L^2(\Omega_1)}$ utilizing the inf-sup condition of the finite element space. For any given $v_h \in X_h(\Omega_1)$, we choose $\xi_h = \text{Ext}(v_h |_{\Gamma})$ and $q_h=0$ in \eqref{weak_Ritz}, with $\text{Ext}(v_h |_{\Gamma})\in X_h(\Omega_2)$ being an extension of $v_h|_\Gamma$ to domain $\Omega_2$ satisfying the following estimate: 
\begin{equation}\label{H1-Ext}
    \begin{aligned}
        \|\xi_h\|_{H^1(\Omega_2)} \le C\|v_h\|_{H^1(\Omega_1)}.
    \end{aligned}
\end{equation}
Such an extension exists. For example, we can extend $v_h |_{\Gamma}$ to $\Omega_2$ using any bounded operator from $H^{1/2}(\Gamma)$ to $H^1(\Omega_2)$ and then project the extended function onto $X_h(\Omega_2)$ using any projection operator which is bounded in the $H^1$ norm. 
Then, for a given time $t\geq 0$, we obtain the following relation: 
\begin{equation}
    \begin{aligned}
        (R_h p - I_h p, \nabla \cdot v_h)_{\Omega_2} 
        & = (\nabla (R_h \eta - I_h \eta), \nabla \xi_h)_{\Omega_2} + (\nabla(I_h \eta - \eta), \nabla \xi_h)_{\Omega_2}  \\
        & \quad \,  + (D(R_h u - I_h u ), D v_h)_{\Omega_1} + (D(I_h u - u), D v_h)_{\Omega_1} \\
        & \quad\, - (I_h p - p, \nabla \cdot v_h)_{\Omega_1}  + (R_ h u - u ,v_h)_{\Omega_1} ,
    \end{aligned}
\end{equation}
where we have omitted the dependence on time $t$ to avoid overloading the notation. 
Employing inequality \eqref{H1-Ext}, we have
\begin{equation}\label{p_L2_infty}
    \begin{aligned}
        \|R_h p -I_h p\|_{L^2(\Omega_1)} &\leq C \sup_{\substack{v_h \in X_h(\Omega_1)\\ 
        v_h\neq 0}}\frac{|(R_h p - I_h p, \nabla \cdot v_h)_{\Omega_2}|}{\|v_h\|_{H^1(\Omega_1)}}\\
        & \leq C\|\nabla (R_h \eta - I_h \eta)\|_{L^2(\Omega_2)} + C\|R_ h u- u\|_{H^1(\Omega_1)} + Ch^k ,
    \end{aligned}
\end{equation}
which implies the following result:  
\begin{equation}\label{p_L2}
    \begin{aligned}
        \|R_h p -I_h p\|_{L^2(0,t; L^2(\Omega_1))} \leq C \|\nabla (R_h \eta - I_h \eta)\|_{L^2(0,t; L^2(\Omega_2))} + C\|R_ h u- u\|_{L^2(0,t; H^1(\Omega_1))} + Ch^k . 
    \end{aligned}
\end{equation}
Combining the estimates in \eqref{L2H1-left}--\eqref{p_L2}, we have
\begin{equation}\label{L2_H1}
    \begin{aligned}
        &\sup_{t\in [0,T]} \|\nabla (R_h \eta - I_h \eta)(t)\|_{L^2(\Omega_1)}^2 + \|R_h u - I_h u\|_{L^2(0,T; H^1(\Omega_1))}^2 \\
        &\leq C \|\nabla (R_h \eta - I_h \eta)(0)\|_{L^2(\Omega_1)}^2 + Ch^{2k} . 
    \end{aligned}
\end{equation}
In order to estimate $\|\nabla (R_h \eta - I_h \eta)(0)\|_{L^2(\Omega_1)}$, we proceed by considering the elliptic problem described in equation \eqref{eta_indep}. This includes adhering to the interface condition where $R_h\eta(0) = I_h \eta(0)$ on $\Gamma$. Therefore, the following equation holds for $\xi_h \in X_h(\Omega_2)$: 
\begin{equation}\label{init-H1}
    \begin{aligned}
        \left\{
        \begin{aligned}
             (\nabla (R_h \eta - I_h \eta)(0), \nabla \xi_h)_{\Omega_2} &= (\nabla (\eta - I_h \eta)(0), \nabla \xi_h)_{\Omega_2}&&\mbox{in}~~ \Omega_2\\
             (R_h \eta - I_h \eta)(0) &= 0 &&\mbox{on}~~ \partial\Omega_2\\
        \end{aligned}
        \right.
    \end{aligned}
\end{equation}
Choosing $\xi_h = (R_h \eta - I_h \eta)(0)$ in \eqref{init-H1} allows us to derive the following inequality:
\begin{equation}\label{init-H1-eta}
    \|\nabla (R_h \eta - I_h \eta)(0)\|_{L^2(\Omega_2)} \leq Ch^k.
\end{equation}
By incorporating the estimate of $\|\nabla (R_h \eta - I_h \eta)(0)\|_{L^2(\Omega_2)}$ from \eqref{init-H1-eta} into inequality \eqref{L2_H1}, we obtain 
\begin{equation}\label{L2_H1_inequality}
    \begin{aligned}
        \sup_{t\in [0,T]} \|\nabla (R_h \eta - I_h \eta)(t)\|_{L^2(\Omega_2)} + \|R_h u - I_h u\|_{L^2(0,T;H^1(\Omega_1))} & \leq Ch^k , \\
        \sup_{t\in [0,T]} \|\nabla (R_h \eta -  \eta)(t)\|_{L^2(\Omega_2)} + \|R_h u - u\|_{L^2(0,T;H^1(\Omega_1))} & \leq Ch^k . 
    \end{aligned}
\end{equation}


\subsection{Estimate of $\|R_h u - u\|_{L^\infty(0,T;H^1(\Omega_1))}$}

We have obtained estimates for $\|\nabla (R_h \eta - \eta)\|_{L^\infty(0,T;L^2(\Omega_2))}$ and $\|R_h u - u\|_{L^2(0,T;H^1(\Omega_1))}$. Our next objective is to improve the estimate of $R_h u - u$ from the $L^2(0,T;H^1(\Omega_1)^d)$ norm to the $L^\infty(0,T;H^1(\Omega_1)^d)$ norm. This improvement can be achieved by extending the test function $v_h$ from $\Omega_1$ to $\Omega_2$ in the definition of the dynamic Ritz projection in \eqref{weak_Ritz}. Namely, we choose test functions $v_h = R_h u - I_h u$, $q_h = R_h p - I_h p$ and $\xi_h = \text{Ext}(v_h|_{\Gamma})$, which satisfies the condition $\|\xi_h\|_{H^1(\Omega_2)} \le C\|v_h\|_{H^1(\Omega_1)}$. Then, by utilizing the bound of $\|R_ h p - I_h p\|_{L^2(\Omega_1)}$ in \eqref{p_L2_infty}, we obtain
\begin{align}
        & \quad \,\|D(R_h u - I_h u)\|_{L^2(\Omega_1)}^2 + \|R_ h u - I_h u \|_{L^2(\Omega_1)}^2 \notag\\
        &= (\nabla ( \eta - R_h \eta), \nabla \xi_h)_{\Omega_2} + (R_h p - p, \nabla \cdot(R_h u - I_h u))_{\Omega_1} \notag\\
        & \quad \, + (\nabla \cdot (u - R_h u ), R_h p - I_h p)_{\Omega_1}+ (D(u - I_h u), D(R_h u - I_h u))_{\Omega_1} \notag\\
        & \quad\, + (u-I_h u, R_h u - I_h u)_{\Omega_1} \notag\\
        & \leq Ch^k (\|\nabla (R_h \eta - I_h \eta)\|_{L^2(\Omega_2)}) + \|R_h u -I_h u\|_{H^1(\Omega_1)} + Ch^k\|R_h p - I_h p\|_{L^2(\Omega_1)} + Ch^{2k} \notag\\
        & \leq Ch^k (\|\nabla (R_h \eta - I_h \eta)\|_{L^2(\Omega_2)} + \|R_ h u - I_h u \|_{H^1(\Omega_1)} + h^k) \notag\\
        & \leq Ch^{2k} + Ch^k \|R_h u - I_h u\|_{H^1(\Omega_1)}.
\end{align}
The last inequality uses the estimate of $\|\nabla (R_h \eta - I_h \eta)\|_{L^\infty(0,T;L^2(\Omega_2))} $ obtained in \eqref{L2_H1_inequality}. Then, by applying Korn's inequality, which controls $\|R_ h u - I_h u \|_{H^1(\Omega_1)}$ by 
$$ \|D(R_h u - I_h u)\|_{L^2(\Omega_1)} + \|R_ h u - I_h u \|_{L^2(\Omega_1)} ,$$ we obtain the following estimate:
\begin{equation}\label{u_h1_Linfty}
    \begin{aligned}
        \sup_{t\in[0,T]}( \|(R_h u - I_h u)(t)\|_{H^1(\Omega_1)} + \|(R_h u - u)(t)\|_{H^1(\Omega_1)}  )&\leq C h^k . 
    \end{aligned}
\end{equation}
This, together with \eqref{p_L2_infty}, implies that 
\begin{equation}\label{p_Lifnty_L2}
    \begin{aligned}
        \sup_{t\in[0,T]} (\|R_h p - I_h p\|_{L^2(\Omega_1)} +  \|R_h p - p\|_{L^2(\Omega_1)}) \le Ch^k.
    \end{aligned}
\end{equation}

\subsection{Estimates of $\|\nabla \partial _t ( R_h \eta -  \eta)\|_{L^\infty(0,T;L^2(\Omega_2))}$ and $\|\partial _t ( R_h u -  u)\|_{L^\infty(0,T;H^1(\Omega_1))}$}

We differentiate equation \eqref{weak_Ritz} with respect to time $t$ to derive the following new equation:
\begin{equation}\label{Ritz_diff}
    \begin{aligned}
        (\nabla \partial_t(R_h \eta - \eta), \nabla \xi_h)_{\Omega_2} &+ (D(\partial_t(R_h u - u)), D(v_h))_{\Omega_1} + (\partial_t(R_h u - u), v_h)_{\Omega_1 }\\
        & - (\partial_t (R_h p - p), \nabla \cdot v_h )_{\Omega_1} + (\nabla \cdot \partial_t(R_h u - u ), q_h)_{\Omega_1}= 0
    \end{aligned}
\end{equation}
For a given time $t\geq 0$, choosing $v_h = 0 $ and $ q_h=0$ in $\Bar{\Omega}_1$, and choosing $\xi_h$ satisfying $\xi_h = 0 $ on $\Gamma$, we have
\begin{equation}\label{partial_eta_indep}
    \left\{
    \begin{aligned}
    &(\nabla \partial_t (R_h \eta - \eta), \nabla \xi_h)_{\Omega_2} = 0 &&\forall\, \xi_h \in X_h(\Omega_2)\,\,\,\mbox{such that}\,\,\,\xi_h=0\,\,\,\mbox{on}\,\,\,\Gamma, \\
    & \partial_t (R_h \eta - \eta) = R_h u - u &&\mbox{ on}\,\,\, \Gamma .
    \end{aligned}
    \right.
\end{equation}
Let us denote by ${\rm Ext}(R_h u - I_hu)\in X_h(\Omega_2)$ an $H^1$-bounded extension of $R_h u - I_hu$ from $\Omega_1$ to $\Omega_2$. Then $\xi_h= \partial_t (R_h \eta -I_h\eta) - {\rm Ext}(R_h u - I_hu)$ satisfies $\xi_h=0$ on $\Gamma$. Choosing this $\xi_h$ in \eqref{partial_eta_indep} leads to the following result: 
\begin{equation}\label{partial_eta_indep_2}
\sup_{t\in[0,T]} \|\nabla \partial_t (R_h \eta - \eta)(t)\|_{L^2(\Omega_2)} 
\le C \sup_{t\in[0,T]}  \|(R_h u - I_hu)(t)\|_{H^1(\Omega_1)} \le C h^k,
\end{equation}
where the last inequality uses \eqref{u_h1_Linfty}. 

Since \eqref{Ritz_diff} only differs from \eqref{weak_Ritz} by an additional time derivative, by choosing test functions $\xi_h = \partial_{tt} (R_h \eta - I_h \eta)$, $v_h = \partial_t(R_h u - I_h u)$, and $q_h = \partial_t(R_h p - I_h p)$ that satisfy the interface condition $\xi_h = v_h$ on $\Gamma$ in the weak formulation \eqref{Ritz_diff}, and utilizing a methodology akin to the derivations in \eqref{L2H1-eqn}--\eqref{L2_H1}, we can obtain the same estimate as \eqref{L2_H1} but with an additional time derivative on the left-hand side, i.e., 
%
\begin{equation}\label{diff_H1_eta}
    \begin{aligned}
        &\sup_{t\in[0,T]}\|\nabla \partial_t (R_h \eta - \eta)(t)\|_{L^2(\Omega_2)}^2 + \| \partial_t (R_h u - u)\|_{L^2 (0,T;H^1(\Omega_1))}^2 \\
        &\le C \|\nabla \partial_t (R_h \eta - \eta)(0)\|_{L^2(\Omega_2)}^2  + C h^{2k} \\
        &\le C h^{2k} , 
    \end{aligned}
\end{equation}
where the last inequality follows from \eqref{partial_eta_indep_2}. 

Furthermore, by selecting test functions $v_h = \partial_t (R_h u - I_h u)$, $q_h = \partial_t (R_h p - I_h p)$ and $\xi_h = \text{Ext}(v_h|_{\Gamma})$ in the weak formulation \eqref{Ritz_diff}, where $\xi_h $ is an extension of $v_h|_{\Gamma}$ to $\Omega_2$ satisfying the condition $\|\xi_h\|_{H^1(\Omega_2)} \le C\|v_h|_\Gamma \|_{H^{1/2}(\Gamma)} \le C\|v_h\|_{H^1(\Omega_1)}$, we can derive the following improved estimate (the pressure term needs to be estimated by using the inf-sup condition):
\begin{align}\label{diff_H1_u}
        &\sup_{t\in[0,T]}
        ( \|\partial_t (R_h u - u)(t)\|_{ H^1(\Omega_1)} + \|\partial_t (R_h p - p)(t)\|_{ L^2(\Omega_1)} ) \notag\\ 
        &\le C \sup_{t\in[0,T]}\|\nabla \partial_t (R_h \eta - \eta)(t)\|_{L^2(\Omega_2)} \le C h^k ,
\end{align}
where the last inequality follows from \eqref{diff_H1_eta}. 

\subsection{Estimates of $\|\nabla \partial _{t}^j ( R_h \eta -  \eta)\|_{L^\infty(0,T;L^2(\Omega_2))}$ and $\| \partial _{t}^j( R_h u -  u)\|_{L^\infty(0,T;H^1(\Omega_1))}$ for $j\geq 1$}

Similarly, we differentiate equation \eqref{Ritz_diff} with respect to $t$ to derive the following new equation:
\begin{equation}\label{Ritz_diff_twice}
    \begin{aligned}
        (\nabla \partial_{tt}(R_h \eta - \eta), \nabla \xi_h)_{\Omega_2} &+ (D(\partial_{tt}(R_h u - u)), D(v_h))_{\Omega_1} + (\partial_{tt}(R_h u - u), v_h)_{\Omega_1 }\\
        & - (\partial_{tt} (R_h p - p), \nabla \cdot v_h )_{\Omega_1} + (\nabla \cdot \partial_{tt}(R_h u - u ), q_h)_{\Omega_1}= 0
    \end{aligned}
\end{equation}
In the same way as we derive \eqref{partial_eta_indep_2}--\eqref{diff_H1_u} from \eqref{Ritz_diff}, by utilizing the estimates in \eqref{partial_eta_indep_2}--\eqref{diff_H1_u} we can derive the following result from \eqref{Ritz_diff_twice}: 
\begin{align}\label{diff_H1_eta_twice} 
\sup_{t\in[0,T]}\|\nabla \partial_{tt} (R_h \eta - \eta)(t)\|_{L^2(\Omega_2)} + \sup_{t\in[0,T]}\|\partial_{tt} (R_h u - u)(t)\|_{ H^1(\Omega_1) } &\le C h^k . 
\end{align}

Then, we differentiate equation \eqref{Ritz_diff_twice} with respect to $t$ to derive the following new equation:
\begin{align}\label{Ritz_diff_triple}
    \begin{aligned}
        (\nabla \partial_{ttt}(R_h \eta - \eta), \nabla \xi_h)_{\Omega_2} &+ (D(\partial_{ttt}(R_h u - u)), D(v_h))_{\Omega_1} + (\partial_{ttt}(R_h u - u), v_h)_{\Omega_1 }\\
        & - (\partial_{ttt} (R_h p - p), \nabla \cdot v_h )_{\Omega_1} + (\nabla \cdot \partial_{ttt}(R_h u - u ), q_h)_{\Omega_1}= 0
    \end{aligned}
\end{align}
In the same way as we derive \eqref{partial_eta_indep_2}--\eqref{diff_H1_u} from \eqref{Ritz_diff}, by utilizing the estimates in \eqref{diff_H1_eta_twice} we can derive the following result from \eqref{Ritz_diff_triple}: 
\begin{align}\label{diff_H1_u_triple}
\sup_{t\in[0,T]}\|\nabla \partial_{ttt} (R_h \eta - \eta)(t)\|_{L^2(\Omega_2)} + \sup_{t\in[0,T]}\|\partial_{ttt} (R_h u - u)(t)\|_{ H^1(\Omega_1) } &\le C h^k . 
\end{align}
This proves the first result of Theorem \ref{THM-Ritz}. The proof provided here is for $j=2,3$, but the bound can be extended to all $j$ using induction.

\subsection{Estimates of $\|(R_h \eta -  \eta)(0)\|_{L^2(\Omega_2)}$ and $\|(R_h u -  u)(0)\|_{L^2(\Omega_1)}$} \label{section:L2-at-t=0}

Note that the following initial condition is given on $\Gamma$: $R_h\eta=I_h\eta$ on $\Gamma$. This allows us to establish an optimal-order estimate of $\|(R_h \eta -  \eta)(0)\|_{L^2(\Omega_2)}$ by considering the following simple dual problem: 
\begin{equation}\label{dual-at-t=0}
    \left\{
    \begin{aligned}
         - \Delta \bar w &= (R_h \eta - \eta)(0) &&{\rm in}~~ \Omega_2\\
         \partial_n\bar w &= 0 &&{\rm on}~~ \partial\Omega_2 \backslash\Gamma \\
         \bar w &= 0 &&{\rm on}~~ \Gamma ,
    \end{aligned}
    \right.
\end{equation}
where the boundary condition guarantees that $I_h\bar w=0$ on $\Gamma$. This and \eqref{eta_indep} imply that  
\begin{align}\label{Ritz-I_h-barw}
(\nabla  I_h \bar w, \nabla (R_h \eta - \eta)(0))_{\Omega_2}=0 . 
\end{align} 
Under Assumption \ref{AS1}, the solution of \eqref{dual-at-t=0} satisfies the following $H^2$ regularity estimate: 
\begin{align}\label{H2-barw}
\|\bar w\|_{H^2(\Omega_2)}\le C\|(R_h \eta - \eta)(0)\|_{L^2(\Omega_2)}.
\end{align} 
By using this result, testing equation \eqref{dual-at-t=0} with $(R_h \eta - \eta)(0)$ and using \eqref{Ritz-I_h-barw}--\eqref{H2-barw} with $R_h \eta=I_h\eta $ on $\Gamma$, we obtain 
\begin{equation*}
    \begin{aligned}
        \|(R_h \eta - \eta)(0)\|_{L^2(\Omega_2)}^2 
        &= (\nabla \bar w, \nabla (R_h \eta - \eta)(0))_{\Omega_2} + (\partial_n\bar w, (R_h \eta - \eta)(0))_\Gamma \\
        &= (\nabla (\bar w - I_h \bar w), \nabla (R_h \eta - \eta)(0))_{\Omega_2} + (\partial_n\bar w, (I_h \eta - \eta)(0))_\Gamma \\
        & \le Ch^k h \|\bar w\|_{H^2(\Omega_2)} \|\eta\|_{H^{k+1}(\Omega_2)} + Ch^{k+1} \|\partial_n\bar w\|_{L^2(\Gamma)} \|\eta\|_{H^{k+1}(\Gamma)} \\
        & \le Ch^{k+1} \|\bar w\|_{H^2(\Omega_2)} \\
        & \le Ch^{k+1}\|(R_h \eta - \eta)(0)\|_{L^2(\Omega_2)}, 
    \end{aligned}
\end{equation*}
for $\eta$ with sufficient smoothness.
This proves that 
\begin{equation}\label{Rheta-eta-at-t=0}
        \|(R_h \eta - \eta)(0)\|_{L^2(\Omega_2)} \le Ch^{k+1} . 
\end{equation}

An optimal-order estimate of $\|(R_h u -  u)(0)\|_{L^2(\Omega_1)}$ can be obtained by considering the following dual problem: 
\begin{equation}\label{dual-u-l2}
    \left\{
    \begin{aligned}
         -\nabla \cdot (D\bar v + \bar qI) + \bar v &= R_h u - u &&{\rm in}~~ \Omega_1\\
         \nabla \cdot \bar v &= 0 &&{\rm in}~~ \Omega_1\\
        (D\bar v + \bar qI)n &= 0 &&{\rm on}~~ \partial\Omega_1 \supset \Gamma . 
    \end{aligned}
    \right.
\end{equation}
According to  Assumption \ref{AS2} on the $H^2$ regularity estimate for the Stokes equations, the solution of \eqref{dual-u-l2} satisfies the following estimate:
\begin{equation}\label{dual-u-l2-reg}
    \begin{aligned}
        \|\bar v\|_{H^2(\Omega_1)} + \|\bar q\|_{H^1(\Omega_1)} &\le C\|R_h u - u\|_{L^2(\Omega_1)} . 
    \end{aligned}
\end{equation}
Testing the first equation of \eqref{dual-u-l2} with $R_h u - u $ and applying integration by parts, and subtracting $I_h\bar v$ and $P_h \bar q$ from $\bar v$ and $\bar q$, respectively, in the resulting weak formulation by using relation \eqref{weak_Ritz}, we obtain 
\begin{equation}\label{Rhu-u-L2}
    \begin{aligned}
        \|R_h u - u\|_{L^2(\Omega_1)}^2  &= (D(\bar v - I_h \bar v), D(R_h u - u))_{\Omega_1} + (\nabla \cdot (R_h u - u), \bar q - P_h \bar q)_{\Omega_1}\\
        & \quad \, + (\bar v - I_h \bar v, R_h u - u)_{\Omega_1} + (\nabla \cdot (\bar v- I_h \bar v), p - R_h p)_{\Omega_1}\\
        & \quad \, - (\nabla (R_h \eta - \eta), \nabla (I_h \bar w - \bar w))_{\Omega_2} - (\nabla (R_h \eta - \eta), \nabla \bar w)_{\Omega_2}\\
        & =: J_1 + J_2 + J_3 + J_4 + J_5 + J_6 ,
    \end{aligned}
\end{equation}
where $\bar w$ can be any $H^1$ function on $\Omega_1$ satisfying $\bar w=v$ on $\Gamma$. We simply choose $\bar w = \text{Ext}(\bar v | _{\Gamma})$ to be an extension of $\bar v|_\Gamma$ satisfying interface condition $\partial_n \bar w= 0$ on $\partial\Omega_2\supset\Gamma$, as well as the following $H^2$ estimate (such an extension exists under Assumption \ref{AS-Ext}): 
$$
\|\bar w\|_{H^2(\Omega_2)} \le C\|\bar v\|_{H^{\frac{3}{2}}(\Gamma)} \le C\|\bar v\|_{H^2(\Omega_1)} . 
$$
Then, from the $H^1$ error estimates in \eqref{L2_H1_inequality} and \eqref{u_h1_Linfty}, and the regularity estimate of $\bar v$ and $\bar q$ in \eqref{dual-u-l2-reg}, we obtain the following estimates:
\begin{equation}
    \begin{aligned}
        |J_1 + J_2 + J_3 + J_4 + J_5 | \le Ch^k h (\|\bar v\|_{H^2(\Omega_1)} + \|\bar q\|_{H^1(\Omega_1)} ) \le Ch^{k+1}\|R_h u  - u\|_{L^2(\Omega_1)} . 
    \end{aligned}
\end{equation}
The last term in \eqref{Rhu-u-L2} can be estimated by using integration by parts and the boundary condition $\partial_nw=0$ on $\partial\Omega_2$ (thus no boundary term is generated when using integration by parts), i.e., 
\begin{equation}
    \begin{aligned}
        |J_6| = |(R_h \eta - \eta, \Delta \bar w)_{\Omega_2}| 
        &\le C\|R_h \eta - \eta\|_{L^2(\Omega_2)}\|\bar w\|_{H^2(\Omega_2)} \\
        &\le C\|R_h \eta - \eta\|_{L^2(\Omega_2)}\|R_h u  - u\|_{L^2(\Omega_1)}.
    \end{aligned}
\end{equation}
Substituting the estimates of $J_m$, $m=1,\dots,6$, into \eqref{Rhu-u-L2}, we obtain 
\begin{equation}\label{u-to-eta}
    \begin{aligned}
        \|R_h u - u\|_{L^2(\Omega_1)} \le Ch^{k+1} + C\|R_h \eta - \eta\|_{L^2(\Omega_2)} 
        \quad\mbox{at any $t\in[0,T]$} . 
    \end{aligned}
\end{equation}
This result, together with \eqref{Rheta-eta-at-t=0}, implies that 
\begin{equation}\label{Rhu-u-at-t=0}
        \|(R_h u - u)(0)\|_{L^2(\Omega_1)} \le Ch^{k+1} . 
\end{equation}

\subsection{Space-time dual problem of the dynamic Ritz projection}

In order to obtain optimal-order error estimates of the dynamic Ritz projection in the $L^\infty(0,T;L^2)$ norm, we consider the following space-time dual problem (backward in time) in the space-time domain $\Omega\times[0,t_*]$ with an arbitrary fixed $t_*\in(0,T]$, with initial condition at $t=t_*$: 
\begin{equation}\label{dual_t}
    \left\{
    \begin{aligned}
         - \Delta w &= \partial_t (\eta - R_h \eta)&&{\rm in}~~ \Omega_2\\
        \nabla \cdot (Dv+qI) - v & = \partial_t (u - R_h u) &&{\rm in}~~ \Omega_1\\
        \nabla \cdot v & = 0 &&{\rm in}~~ \Omega_1
    \end{aligned}
    \right.
\end{equation}
with interface conditions
\begin{equation}\label{dual_t_interface}
    \left\{
    \begin{aligned}
        \partial_t w &= v &&{\rm on}~~ \Gamma\\
        w(\cdot, t_*) &= 0 &&{\rm on}~~ \Gamma\\
        - \partial_nw &= (Dv + qI) n &&{\rm on}~~ \Gamma
    \end{aligned}
    \right.
\end{equation}
and boundary conditions 
\begin{equation}\label{dual_t_boundary}
    \left\{
    \begin{aligned}
        \partial_nw &= 0 &&{\rm on}~~ \partial \Omega_2 \backslash \Gamma\\
        (Dv+qI)n &= 0 &&{\rm on}~~ \partial \Omega_1 \backslash \Gamma . 
    \end{aligned}
    \right.
\end{equation}


We estimate $\|v\|_{H^2(\Omega_1)} + \|q\|_{H^1(\Omega_1)} $ and $\|w\|_{H^2(\Omega_2)} $ as follows, by using the $H^2$ regularity estimates in \eqref{H2-v} and \eqref{H2-w}, i.e., 
\begin{align}\label{v_H2}
        &\|v\|_{H^2(\Omega_1)} + \|q\|_{H^1(\Omega_1)} \notag\\
        &\le C\|\partial_t(R_h u - u)\|_{L^2(\Omega_1)} + C\|(D(v)+qI)n\|_{H^{\frac{1}{2}}(\Gamma)} 
        \notag\\
        &= C\|\partial_t(R_h u - u)\|_{L^2(\Omega_1)} + C\|\partial_nw\|_{H^{\frac{1}{2}}(\Gamma)} 
        \quad\,\mbox{(since $-\partial_nw=(D(v)+qI)n$ on $\Gamma$)}\notag\\
        &\le C\|\partial_t(R_h u - u)\|_{L^2(\Omega_1)} + C\| w\|_{H^{2}(\Omega_2)} 
         \qquad\mbox{(trace inequality)}\notag\\
        & \le C\|\partial_t(R_h u - u)\|_{L^2(\Omega_1)} +C \|\partial_t(R_h \eta - \eta)\|_{L^2(\Omega_2)}+ C\|w\|_{H^{\frac{3}{2}}(\Gamma)} 
        \quad\mbox{(here \eqref{H2-w} is used)} \notag\\
        &\le C\|\partial_t(R_h u - u)\|_{L^2(\Omega_1)} + C \|\partial_t(R_h \eta - \eta)\|_{L^2(\Omega_2)}
            + C\Big\|\int_t^{t_*} v(s)\d s \Big\|_{H^{\frac{3}{2}}(\Gamma)} \notag\\
        &\hspace{15pt}\mbox{(note that $w(t)=\int_t^{t_*} v(s)\d s$ as a result of $ \partial_t w=v$ and $w(\cdot,t_*)=0$ on $\Gamma$)}\notag\\
        &\le C\|\partial_t(R_h u - u)\|_{L^2(\Omega_1)} + C\|\partial_t(R_h \eta - \eta)\|_{L^2(\Omega_2)}+ C\int_t^{t_*} \|v(s)\|_{H^2(\Omega_1)} \d s, 
\end{align}
which, together with Gronwall's inequality with initial condition at $t=t_*$, implies that 
\begin{equation}\label{v_LinftyH2}
    \begin{aligned}
        &\|v\|_{L^\infty(0,t_*; H^2(\Omega_1))} +  \|q\|_{L^\infty(0,t_*; H^1(\Omega_1))} \\
        &\le C ( \|\partial_t(R_h \eta - \eta)\|_{L^\infty(0,t_*;L^2(\Omega_2))} + \|\partial_t(R_h u - u)\|_{L^\infty(0,t_*;L^2(\Omega_1))} ) .
    \end{aligned}
\end{equation}
Additionally, the regularity estimates for $w$ at any given time $t\in[0,T]$ can be derived from the regularity estimates in \eqref{H2-w}, i.e., 
\begin{align*}
    \|w\|_{H^2(\Omega_2)} 
    &\le C\|\partial_t(R_h \eta - \eta)\|_{L^2(\Omega_2)} + C\|\partial_n w\|_{H^{\frac{1}{2}}(\Gamma)} \notag\\
    &= C\|\partial_t(R_h \eta - \eta)\|_{L^2(\Omega_2)} + C\|(D(v)+qI)n\|_{H^{\frac{1}{2}}(\Gamma)} \notag\\
    &\le C\|\partial_t(R_h \eta - \eta)\|_{L^2(\Omega_2)} + C\|v\|_{H^2(\Omega_1)} + C\|q\|_{H^1(\Omega_1)} ,
\end{align*}
which implies that 
\begin{equation}\label{w_H2}
        \|w\|_{L^\infty(0,t_*;H^2(\Omega_2))}   \le C( \|\partial_t(R_h \eta - \eta)\|_{L^\infty(0,t_*; L^2(\Omega_2))} + \|\partial_t(R_h u - u)\|_{L^\infty(0,t_*;L^2(\Omega_1))} ) .
\end{equation}

\subsection{Estimates of $\|\partial_t ( R_h \eta -  \eta)\|_{L^\infty(0,T;L^2(\Omega_2))}$ and $\|\partial_t (R_h u - u)\|_{L^\infty(0,T;L^2(\Omega_1))}$}

We establish the estimates of $\|\partial_t ( R_h \eta -  \eta)\|_{L^\infty(0,T;L^2(\Omega_2))}$ and $\|\partial_t (R_h u - u)\|_{L^\infty(0,T;L^2(\Omega_1))}$ via a duality argument based on the space-time dual problem introduced in \eqref{dual_t}--\eqref{dual_t_boundary}. The weak formulation of \eqref{dual_t}--\eqref{dual_t_boundary} can be written as follows: 
\begin{equation}\label{dual_partial_t}
    \begin{aligned}
        (\partial_t (\eta - R_h \eta), \xi)_{\Omega_2} + (\partial_t (u-R_h u) , \zeta)_{\Omega_1} &= (\nabla w, \nabla \xi)_{\Omega_2} - (Dv , D\zeta)_{\Omega_1} - (q, \nabla \cdot \zeta)_{\Omega_1} \\
        & \quad\,- (v, \zeta)_{\Omega_1} + (\nabla \cdot v, \psi)_{\Omega_1} ,
    \end{aligned}
\end{equation}
which holds for test functions satisfying interface condition $\xi = \zeta$ on $\Gamma$. 
Choosing test functions $\xi = \partial_{tt} (\eta - R_h \eta)$, $\zeta = \partial_t (u - R_h u)$ and $\psi = \partial_t (p - R_h p)$ which satisfy the interface condition $\xi = \zeta$ on $\Gamma$ in \eqref{dual_partial_t}, we obtain
\begin{equation}\label{dual_identity_partial_t}
    \begin{aligned}
        & \quad \,\frac{1}{2}\frac{\d}{\d t} \|\partial_t(R_h \eta - \eta)\|_{L^2(\Omega_2)}^2 + \|\partial_t(R_h u - u)\|_{L^2(\Omega_1)}^2 \\
        &= - (\nabla w, \nabla \partial_{tt} (R_h \eta - \eta))_{\Omega_2} + (Dv , D(\partial_t(R_h u - u)))_{\Omega_1} \\
        & \quad \,+ (q, \nabla \cdot \partial_t(R_h u - u))_{\Omega_1} +(v, \partial_t(R_h u - u))_{\Omega_1}\\
        & \quad \, + (\nabla \cdot v, \partial_t(p - R_h p))_{\Omega_1}\\
        & =: K_1 + K_2 + K_3 + K_4 + K_5 . 
    \end{aligned}
\end{equation}
Integrating \eqref{dual_identity_partial_t} from $0$ to $t_*$ with respect to time, we obtain 
\begin{equation}\label{dual_identity_partial_t_integrate}
    \begin{aligned}
        & \quad \, \frac{1}{2}\|\partial_t (R_h \eta - \eta)(t_*)\|_{L^2(\Omega_2)}^2 - \frac{1}{2}\|\partial_t (R_h \eta - \eta)(0)\|_{L^2(\Omega_2)}^2 + \int_0^{t_*} \|\partial_t(R_h u - u)\|_{L^2(\Omega_1)}^2 \d t\\
        & = \int_0^{t_*} (K_1 + K_2 + K_3 + K_4 + K_5)\d t , 
    \end{aligned}
\end{equation}
where 
\begin{equation}\label{def-K1-K5}
    \begin{aligned}
        & K_1 = -\, \partial_t (\nabla w, \nabla \partial_t(R_h \eta - \eta))_{\Omega_2} + (\nabla \partial_t (w - I_h w), \nabla \partial_t (R_h \eta - \eta))_{\Omega_2} \\
        & \qquad\,\,\, + (\nabla \partial_t I_h w, \nabla \partial_t (R_h \eta - \eta))_{\Omega_2}
        =: K_{11} + K_{12} + K_{13}\\
        & K_2 = (D(v - I_h v), D(\partial_t(R_h u - u)))_{\Omega_1} + (D I_h v, D(\partial_t(R_h u - u)))_{\Omega_1} =: K_{21} + K_{22}\\
        & K_3 = (q - P_h q, \nabla \cdot \partial_t(R_h u - u))_{\Omega_1} + (P_h q, \nabla \cdot \partial_t(R_h u - u))_{\Omega_1} =: K_{31} + K_{32}\\
        & K_4 = (v - I_h v, \partial_t(R_h u - u))_{\Omega_1} + (I_h v, \partial_t(R_h u - u))_{\Omega_1} =: K_{41} + K_{42}\\
        & K_5 = -(\nabla \cdot(v - I_h v), \partial_t(R_hp -  p))_{\Omega_1} - (\nabla \cdot I_h v, \partial_t(R_hp -  p))_{\Omega_1} =: K_{51} + K_{52}.
    \end{aligned}
\end{equation}
Since $\partial_tw=v$ on $\Gamma$, it follows that $\partial_tI_hw=I_hv$ on $\Gamma$. Under this condition, from the weak formulation in \eqref{Ritz_diff} we can derive that 
\begin{equation}\label{K13=0}
    \begin{aligned}
        K_{13} + K_{22} + K_{32} + K_{42} + K_{52} &= 0 .
    \end{aligned}
\end{equation}
By replicating the approach used to establish \eqref{p_L2_infty} and applying the inf-sup conditions to the weak formulation in \eqref{Ritz_diff} (the temporally differentiated weak formulation), we arrive at the following result for any given time $t\in[0,T]$: 
\begin{equation}\label{p_L2_t_infty}
    \|\partial_t (R_h p - p)\|_{L^2(\Omega_1)} \le Ch^k + \|\partial_t (R_h u - u)\|_{H^1(\Omega_1)} + \|\nabla \partial_t (R_h \eta - \eta)\|_{L^2(\Omega_2)} \le Ch^k . 
\end{equation}
By utilizing the regularity estimates \eqref{v_LinftyH2}--\eqref{w_H2}, we obtain
\begin{equation}\label{int-K21}
    \begin{aligned}
        & \int_0^{t_*} |K_{21} + K_{31} + K_{41} + K_{51}| \d t \\
        &\le Ch^{k+1}(\|\partial_t(R_h \eta - \eta)\|_{L^\infty(0,t_*;L^2(\Omega_2))} + \|\partial_t(R_h u - u)\|_{L^\infty(0,t_*;L^2(\Omega_1))}) . 
    \end{aligned}
\end{equation}
By moving $\partial_t$ from $w-I_hw$ to $\partial_t(R_h\eta-\eta)$ in the expression of $K_{12}$, we can rewrite $K_{11} + K_{12}$ as follows: 
\begin{equation*}
    \begin{aligned}
        K_{11} + K_{12} = - \partial_t (\nabla I_h w, \nabla \partial_t(R_h \eta - \eta))_{\Omega_2} - (\nabla (w - I_h w), \nabla \partial_{tt} (R_h \eta - \eta))_{\Omega_2} 
        =: \bar K_{11} + \bar K_{12} . 
    \end{aligned}
\end{equation*}
The term $\bar K_{12}$ can be estimated by using the inequality \eqref{diff_H1_eta_twice}, and the regularity estimates in \eqref{v_LinftyH2}--\eqref{w_H2}, we obtain
\begin{equation}\label{int-K12}
    \begin{aligned}
        \int_0^{t_*} |\bar K_{12}| \d t &\le Ch^{k+1} \|w\|_{L^\infty(0,t_*;H^2(\Omega_2))} \\
        &\le Ch^{k+1} (\|\partial_t(R_h \eta - \eta)\|_{L^\infty(0,t_*;L^2(\Omega_2))} + \|\partial_t(R_h u - u)\|_{L^\infty(0,t_*;L^2(\Omega_1))}) . 
    \end{aligned}
\end{equation}
Furthermore, since $w(t_*)=0$ on $\Gamma$, it follows that \( I_h w(t_*) = 0 \) on $\Gamma$. By choosing $v_h=0$, $q_h=0$ and $\xi_h = I_hw(t_*)$ in the weak formulation \eqref{Ritz_diff} (thus $\xi_h=v_h$ on $\Gamma$), we have 
$$ (\nabla I_h w(t_*), \nabla \partial_t(R_h \eta - \eta)(t_*))_{\Omega_2} = 0 . $$ 
By using this relation, we have 
\begin{equation}\label{K11}
    \begin{aligned}
        \int_0^{t_*} \bar K_{11} \d t &= (\nabla (I_h w(0) - w(0)), \nabla \partial_t (R_h \eta - \eta)(0))_{\Omega_2}  + (\nabla w(0), \nabla \partial_t (R_h \eta - \eta)(0))_{\Omega_2} \\ 
        & =: \bar K_{111} + \bar K_{112} . 
    \end{aligned}
\end{equation}
Employing the $H^1$ error estimates in \eqref{diff_H1_eta}--\eqref{diff_H1_u} and the regularity estimates in \eqref{v_LinftyH2}--\eqref{w_H2}, we have 
\begin{equation}
    \begin{aligned}
        |\bar K_{111}| \le Ch^{k+1} (\|\partial_t(R_h \eta - \eta)\|_{L^\infty(0,t_*;L^2(\Omega_2))} + \|\partial_t(R_h u - u)\|_{L^\infty(0,t_*; L^2(\Omega_1))})
    \end{aligned}
\end{equation}
By using integration by parts, we have
\begin{equation}\label{K112}
    \begin{aligned}
        \bar K_{112} &= - (\Delta w(0), \partial_t (R_h \eta - \eta)(0))_{\Omega_2} - (\partial_n w(0), \partial_t (R_h \eta - \eta)(0))_{\Gamma}\\
        & =  - \|\partial_t (R_h \eta - \eta)(0)\|_{L^2}^2 - (\partial_n w(0), \partial_t (R_h \eta - \eta)(0))_{\Gamma}
    \end{aligned}
\end{equation}
where the first term on the right-hand side of \eqref{K112} can be used to cancel the second term on the left-hand side of  \eqref{dual_identity_partial_t_integrate}. This is key cancellation structure which allows us to establish optimal-order estimates in the $L^\infty(0,T;L^2)$ norm. 

Additionally, using the interface conditions 
$-\partial_nw= (D(v)+qI)n$ and $\partial_t (R_h \eta - \eta)= R_h u - u$ on $\Gamma$, we have 
\begin{equation}\label{K112-2}
    \begin{aligned}
        &\quad \, - (\partial_nw(0), \partial_t (R_h \eta - \eta)(0))_{\Gamma}\\
        & = ((D(v(0)) + q(0)I)n, (R_h u - u)(0))_{\Gamma}\\
        & = (\nabla \cdot (D(v(0)) + q(0)I) - v(0), (R_h u - u)(0))_{\Omega_1} + (v(0), (R_h u - u)(0))_{\Omega_1} \\
        & \quad\, + (D(v(0)), D(R_h u - u)(0))_{\Omega_1} 
                               + (q(0), \nabla \cdot (R_h u - u)(0))_{\Omega_1}\\
        & = -(\partial_t (R_h u - u)(0), (R_h u - u)(0))_{\Omega_1} + ((v - I_h v)(0), (R_h u - u)(0))_{\Omega_1}\\
        & \quad \, + (D(v - I_h v)(0), D(R_h u - u)(0))_{\Omega_1} + ((q - P_h q)(0), \nabla \cdot(R_h u - u)(0))_{\Omega_1}\\
        & \quad \, - (\nabla (R_h \eta - \eta)(0), \nabla \xi_h)_{\Omega_2} 
        + ((R_h p - p)(0), \nabla \cdot I_h v(0) )_{\Omega_2} 
        \quad\mbox{(here \eqref{weak_Ritz_Omega} is used)} \\
        & =: L_1 + L_2 + L_3 + L_4 + L_5 + L_6 , 
    \end{aligned} 
\end{equation}
where $\xi_h = I_h \xi\in X_h(\Omega_2)$ for some $\xi={\rm Ext}(v|_{\Gamma})$ which is an extension of $v|_{\Gamma}$ from $\Gamma$ to $\Omega_2$ such that $\partial_n\xi=0$ on $\partial\Omega_2$ and $\|\xi\|_{H^2(\Omega_2)}\le C\|v|_{\Gamma}\|_{H^{\frac32}(\Gamma)}\le C\|v\|_{H^2(\Omega_1)}$ (such an extension $\xi$ exists under Assumption \ref{AS-Ext}). 

Note that $L_1$ can be estimated by using estimates of $\|(R_h u - u)(0)\|_{L^2(\Omega_1)}$ established in \eqref{Rheta-eta-at-t=0}, which implies that 
\begin{align*}
        L_1 \le Ch^{k+1} \|\partial_t(R_h u - u)\|_{L^\infty(0,t_*;L^2(\Omega_1))} . 
\end{align*}
Since $\nabla\cdot v(0)=0$, it follows that $L_6$ can be rewritten as $ ((R_h p - p)(0), \nabla \cdot[I_h v(0) - v(0)] )_{\Omega_2} $ and then estimated by using \eqref{p_Lifnty_L2} and the regularity estimate of $v$ in \eqref{v_LinftyH2}--\eqref{w_H2}. Additionally, $L_2$, $L_3$ and $L_4$ can be estimated similarly by using the $H^1$ error estimates in \eqref{L2_H1_inequality} and \eqref{u_h1_Linfty} with the regularity estimates of $v$ and $q$ in \eqref{v_LinftyH2}--\eqref{w_H2}, i.e., 
\begin{align*}
        L_2 + L_3 + L_4 + L_6
        &\le Ch^{k+1} (\|v\|_{L^\infty(0,t_*;H^2(\Omega_1))} + \|q\|_{L^\infty(0,t_*;H^1(\Omega_1))}) \\
        &\le Ch^{k+1} (\|\partial_t(R_h \eta - \eta)\|_{L^\infty(0,t_*;L^2(\Omega_2))}
        + \|\partial_t(R_h u - u)\|_{L^\infty(0,t_*;L^2(\Omega_1))}) . 
\end{align*}
Furthermore, $L_5$ can be estimated with integration by parts and interface condition $\partial_n \xi=0$ on $\Gamma$, i.e., 
\begin{align*}
L_5 
&= - (\nabla (R_h \eta - \eta)(0), \nabla \xi)_{\Omega_2}  - (\nabla (R_h \eta - \eta)(0), \nabla (\xi_h-\xi))_{\Omega_2} \\
&= ( (R_h \eta - \eta)(0), \Delta \xi)_{\Omega_2}  - (\nabla (R_h \eta - \eta)(0), \nabla (\xi_h-\xi))_{\Omega_2} \\
&\le C\| (R_h \eta - \eta)(0)\|_{L^2(\Omega_2)} \| \xi\|_{L^\infty(0,t_*;H^2(\Omega_2))} 
 + \| R_h \eta - \eta\|_{L^\infty(0,t_*;H^1(\Omega_2))} Ch \| \xi\|_{L^\infty(0,t_*;H^2(\Omega_2))} \\
&\le Ch^{k+1} (\|\partial_t(R_h \eta - \eta)\|_{L^\infty(0,t_*;L^2(\Omega_2))}
        + \|\partial_t(R_h u - u)\|_{L^\infty(0,t_*;L^2(\Omega_1))}) , 
\end{align*}
where the last inequality uses \eqref{Rheta-eta-at-t=0}, \eqref{L2_H1_inequality} and $\| \xi\|_{H^2(\Omega_2)} \le C\| v\|_{H^2(\Omega_1)} $, together with the  estimate of $\| v\|_{H^2(\Omega_1)} $ in \eqref{v_LinftyH2}. 

Substituting the estimates of $L_j$, $j=1,\cdots,6$, into \eqref{K112-2} and \eqref{K112}, we obtain the desired estimate of $\bar K_{112}$. Then, substituting the estimates of $\bar K_{111}$ and $\bar K_{112}$ into \eqref{K11}, we obtain 
\begin{equation}\label{int-K11}
    \begin{aligned}
        \int_0^{t_*} \bar K_{11} \d t 
        &\le Ch^{k+1} (\|\partial_t(R_h \eta - \eta)\|_{L^\infty(0,t_*;L^2(\Omega_2))} + \|\partial_t(R_h u - u)\|_{L^\infty(0,t_*;L^2(\Omega_1))}) \\
        &\quad\, - \|\partial_t (R_h \eta - \eta)(0)\|_{L^2}^2.  
    \end{aligned}
\end{equation} 
Then, substituting the estimates of \eqref{int-K11}, \eqref{int-K12} and \eqref{int-K21} into \eqref{dual_identity_partial_t_integrate}--\eqref{def-K1-K5} and using relation \eqref{K13=0}, we obtain 
\begin{equation}\label{dt(R_heta-eta)}
    \begin{aligned}
        & \quad \, \frac{1}{2}\|\partial_t (R_h \eta - \eta)(t_*)\|_{L^2(\Omega_2)}^2  + \frac{1}{2}\|\partial_t (R_h \eta - \eta)(0)\|_{L^2(\Omega_2)}^2  +  \|\partial_t(R_h u - u)\|_{L^2(0,t_*;L^2(\Omega_1))}^2 \\
        & \le Ch^{k+1} (\|\partial_t(R_h \eta - \eta)\|_{L^\infty(0,t_*;L^2(\Omega_2))} + \|\partial_t(R_h u - u)\|_{L^\infty(0,t_*;L^2(\Omega_1))}) . 
    \end{aligned}
\end{equation}

It remains to estimate $\|\partial_t(R_h u - u)\|_{L^\infty(0,t_*;L^2(\Omega_1))}$ on the right-hand side of \eqref{dt(R_heta-eta)}. This is achieved by considering the following dual problem for any fixed $t\in[0,T]$ (we omit the dependence on $t$ for the simplicity of notation): 
\begin{equation}\label{dual-u-l2_t}
    \left\{
    \begin{aligned}
        -\nabla \cdot (D\bar v + \bar qI) + \bar v & = \partial_t(R_h u - u) &&{\rm in}\,\,\,  \Omega_1\\
        \nabla \cdot \bar v & = 0 &&{\rm in}\,\,\, \Omega_1\\
        (D\bar v + \bar qI)n & = 0 &&{\rm on}\,\,\,  \partial \Omega_1 .\\
    \end{aligned}
    \right.
\end{equation}
The regularity estimates of Stokes equations in \eqref{H2-v} guarantee that 
\begin{equation}\label{dual-u-l2-reg_t}
    \begin{aligned}
        \|\bar v\|_{H^2(\Omega_1)} + \|\bar q\|_{H^1(\Omega_1)} &\le C\|\partial_t(R_h u - u)\|_{L^2(\Omega_1)} . 
    \end{aligned}
\end{equation}
By testing \eqref{dual-u-l2_t} with $\partial_t(R_h u - u) $ and using equation \eqref{Ritz_diff} (which allows us to subtract $I_h\bar v$ from the weak formulation), we obtain 
\begin{equation}\label{dt-Rhu-u-L2}
    \begin{aligned}
        \|\partial_t(R_h u - u)\|_{L^2(\Omega_1)}^2  &= (D(\bar v - I_h \bar v), D\partial_t(R_h u - u))_{\Omega_1} + (\nabla \cdot \partial_t(R_h u - u), \bar q - P_h \bar q)_{\Omega_1}\\
        & \quad \, + (\bar v - I_h \bar v, \partial_t(R_h u - u))_{\Omega_1} + (\nabla \cdot (\bar v- I_h \bar v), \partial_t(p - R_h p))_{\Omega_1}\\
        & \quad \, - (\nabla \partial_t(R_h \eta - \eta), \nabla (I_h \bar w - \bar w))_{\Omega_2} - (\nabla \partial_t(R_h \eta - \eta), \nabla \bar w)_{\Omega_2}\\
        & := M_1 + M_2 + M_3 + M_4 + M_5 + M_6,
    \end{aligned}
\end{equation}
where $\bar w = \text{Ext}(\bar v | _{\Gamma})$ ia an extension satisfying  $\partial_n \bar w= 0$ on $\partial\Omega_2$ and regularity estimate $\|\bar w\|_{H^2(\Omega_2)} \le C\|\bar v\|_{H^{\frac{3}{2}}(\Gamma)} \le C\|\bar v\|_{H^2(\Omega_1)}$ (such an extension $\bar w$ exists under Assumption \ref{AS-Ext}). From the regularity estimate in \eqref{dual-u-l2-reg_t}, we see that
\begin{equation}
    \begin{aligned}
        |M_1 + M_2 + M_3 + M_4 + M_5 | \le Ch^k h (\|\bar v\|_{H^2(\Omega_1)} + \|\bar q\|_{H^1(\Omega_1)} ) \le Ch^{k+1}\|\partial_t(R_h u  - u)\|_{L^2(\Omega_1)}.
    \end{aligned}
\end{equation}
Moreover, $M_6$ can be estimated with integration by parts using the property $\partial_n \bar w= 0$ on $\partial\Omega_2$, i.e., 
\begin{equation}
    \begin{aligned}
        |M_6| = |(\partial_t(R_h \eta - \eta), \Delta \bar w)_{\Omega_1}| 
        &\le C\|\partial_t(R_h \eta - \eta)\|_{L^2(\Omega_2)}\|\bar w\|_{H^2(\Omega_2)} \\
        &\le C\|\partial_t(R_h \eta - \eta)\|_{L^2(\Omega_2)}\|\bar v\|_{H^2(\Omega_1)} \\
        &\le C\|\partial_t(R_h \eta - \eta)\|_{L^2(\Omega_2)}\|\partial_t(R_h u  - u)\|_{L^2(\Omega_1)} ,
    \end{aligned}
\end{equation}
where \eqref{dual-u-l2-reg_t} is used in the last inequality. 
This implies that, after substituting the estimates of $M_j$, $j=1,\dots,6$, into \eqref{dt-Rhu-u-L2}, 
\begin{equation}\label{u-to-eta_t}
    \begin{aligned}
        \|\partial_t(R_h u - u)\|_{L^\infty(0,T;L^2(\Omega_1))} \le Ch^{k+1} + C\|\partial_t(R_h \eta - \eta)\|_{L^\infty(0,T;L^2(\Omega_2))} . 
    \end{aligned}
\end{equation}
Then, substituting \eqref{u-to-eta_t} into the right-hand side of \eqref{dt(R_heta-eta)}, and taking $L^\infty$ norm with respect to $t_*\in(0,T]$, we obtain 
\begin{equation*}
    \begin{aligned}
        \|\partial_t(R_h \eta - \eta)\|_{L^\infty(0,T;L^2(\Omega_2)}^2 
        &\leq  \epsilon \|\partial_t(R_h \eta - \eta)\|_{L^\infty(0,T;L^2(\Omega_2))}^2 + C_\epsilon h^{2k+2} . 
    \end{aligned}
\end{equation*}
This gives us the following estimate: 
\begin{equation}\label{eta-l2-result_t}
    \begin{aligned}
        \|\partial_t(R_h \eta - \eta)\|_{L^\infty(0,T;L^2(\Omega_2)} \le C h^{k+1} . 
    \end{aligned}
\end{equation}
Then, by utilizing \eqref{u-to-eta_t} again, we obtain
\begin{equation}\label{u-l2-result_t}
    \begin{aligned}
        \|\partial_t(R_h u - u)\|_{L^\infty(0,T;L^2(\Omega_1)}\le Ch^{k+1} . 
    \end{aligned}
\end{equation}
This, together with the estimates of $ \|(R_h \eta - \eta)(0)\|_{L^2(\Omega_2)}$ and $ \|(R_h u - u)(0)\|_{L^2(\Omega_1)}$ in \eqref{Rheta-eta-at-t=0} and \eqref{Rhu-u-at-t=0}, respectively, yields the following result:
\begin{equation}\label{linf-l2-result}
    \begin{aligned}
        \|R_h \eta - \eta\|_{L^\infty(0,T;L^2(\Omega_2))}
        + \|R_h u - u\|_{L^\infty(0,T;L^2(\Omega_1))} \le C h^{k+1} . 
    \end{aligned}
\end{equation}

\subsection{Estimates of $\|\partial_{tt} ( R_h \eta -  \eta)\|_{L^\infty(0,T;L^2(\Omega_2))}$ and $\|\partial_{tt} (R_h u - u)\|_{L^\infty(0,T;L^2(\Omega_1))}$}

Optimal-order estimates of $\|\partial_{tt} ( R_h \eta -  \eta)\|_{L^\infty(0,T;L^2(\Omega_2))}$ and $\|\partial_{tt} (R_h u - u)\|_{L^\infty(0,T;L^2(\Omega_1))}$ can be obtained by considering the following space-time dual problem (backward in time) with initial condition at $t=t_*$ (for an arbitrary fixed $t_*\in(0,T]$): 
\begin{equation}\label{dual_tt}
    \left\{
    \begin{aligned}
         - \Delta w &= \partial_{tt} (\eta - R_h \eta)&&{\rm in}~~ \Omega_2\\
        \nabla \cdot (Dv+qI) - v & = \partial_{tt}  (u - R_h u) &&{\rm in}~~ \Omega_1\\
        \nabla \cdot v & = 0 &&{\rm in}~~ \Omega_1
    \end{aligned}
    \right.
\end{equation}
with interface conditions
\begin{equation}\label{dual_tt_interface}
    \left\{
    \begin{aligned}
        \partial_t w &= v &&{\rm on}~~ \Gamma\\
        w(\cdot, t_*) &= 0 &&{\rm on}~~ \Gamma\\
        - \partial_nw &= (Dv + qI) n &&{\rm on}~~ \Gamma
    \end{aligned}
    \right.
\end{equation}
and boundary conditions 
\begin{equation}\label{dual_tt_boundary}
    \left\{
    \begin{aligned}
        \partial_n w &= 0 &&{\rm on}~~ \partial \Omega_2 \backslash \Gamma\\
        (Dv+qI)n &= 0 &&{\rm on}~~ \partial \Omega_1 \backslash \Gamma . 
    \end{aligned}
    \right.
\end{equation}
In the last section, we have obtained optimal-order estimates of $\|\partial_{t} ( R_h \eta -  \eta)\|_{L^\infty(0,T;L^2(\Omega_2))}$ and $\|\partial_{t} (R_h u - u)\|_{L^\infty(0,T;L^2(\Omega_1))}$ by utilizing the optimal-order estimates of $\|(R_h \eta -  \eta)(0)\|_{L^2(\Omega_2)}$ and $\|(R_h u - u)(0)\|_{L^2(\Omega_1)}$ established in Section \ref{section:L2-at-t=0}. 
Optimal-order estimates of $\|\partial_{tt} ( R_h \eta -  \eta)\|_{L^\infty(0,T;L^2(\Omega_2))}$ and $\|\partial_{tt} (R_h u - u)\|_{L^\infty(0,T;L^2(\Omega_1))}$ can be obtained in the same way  by utilizing the optimal-order estimates of $\|\partial_t(R_h \eta -  \eta)(0)\|_{L^2(\Omega_2)}$ and $\|\partial_t(R_h u - u)(0)\|_{L^2(\Omega_1)}$. The latter follow from  \eqref{eta-l2-result_t}--\eqref{u-l2-result_t} and leads to the estimate: 
\begin{equation}\label{eta-l2-result_tt}
    \begin{aligned}
        \|\partial_{tt}(R_h \eta - \eta)\|_{L^\infty(0,T;L^2(\Omega_2))} 
        +  \|\partial_{tt}(R_h u - u)\|_{L^\infty(0,T;L^2(\Omega_1))} 
        \le C h^{k+1} . 
    \end{aligned} 
\end{equation} 
This proves the second result of Theorem \ref{THM-Ritz}. The proof of Theorem \ref{THM-Ritz} is complete.
\hfill\qed

\section{Proof of Theorem \ref{THM-Error}}\label{section:error}

In this section, we prove the error estimates in Theorem \ref{THM-Error} by utilizing the results we have proved in Theorem \ref{THM-Ritz}. 

\subsection{Error equations}
For the simplicity of notation, we denote $R_h \eta(t_{n})$, $R_h w(t_n)$, $R_h u(t_n)$, $R_h p(t_n)$ by $R_h \eta^{n}$, $R_h w^n$, $R_h u^n$, $R_h p^n$, respectively, and define  
$$R_h w := \partial_t R_h \eta 
\quad\mbox{and}\quad 
R_h p^{n+1/2} : = (R_hp(t_n) + R_h p (t_{n+1}))/2 .
$$
Then, replacing \((\eta_h, w_h, u_h, p_h)\) by \((R_h \eta, R_h w, R_h u, R_h p)\) in the weak formulation \eqref{numerical_weak} leads to the following discrete weak formulation with defect terms: 
\begin{subequations}\label{defects-eqn}
    \begin{align}
        & (R_h \eta^{n+1} - R_h \eta^n)/\tau = (R_h w^{n+1} + R_h w^n)/2 + d_w\quad\mbox{on}\,\,\, \Omega_2, \label{defects-eqn-1}\\[5pt]
        & ((R_h w^{n+1} - R_h w^n)/\tau, \xi_h)_{\Omega_{2}} + (\nabla(R_h \eta^{n+1} + R_h \eta^n)/2, \nabla \xi_h)_{\Omega_{2}} \notag\\
        &\quad\, + ((R_h u^{n+1} - R_h u^n)/\tau, v_h)_{\Omega_{1}} + (D((R_h u^{n+1} + R_h u^n)/2) , D(v_h))_{\Omega_{1}} \notag\\
        &\quad\,- (R_h p^{n+1/2} , \nabla \cdot v_h )_{\Omega_{1}}+(\nabla \cdot (R_h u^{n+1} + R_h u^n)/2, q_h)_{\Omega_{1}} \notag\\
        & = 
         (f(t_{n+1/2}), v_h)_{\Omega_{1}} + (g(t_{n+1/2}), v)_{\partial \Omega_1\backslash\Gamma} 
        + d_\eta(\xi_h) + d_u(v_h). \label{defects-eqn-2}
    \end{align}
\end{subequations}
From the definition of the dynamic Ritz projection \eqref{weak_Ritz} and the weak formulation \eqref{exact-weak} of the exact solutions at time $t_n$ and $t_{n+1}$, we can derive the following expressions for the defect terms: 
\begin{subequations}\label{defect-terms}
\begin{align}
    d_w &= ( R_h \eta^{n+1} -  R_h\eta^n)/\tau - \partial_t R_h \eta(t_{n+1/2}) \notag\\
    &\quad\, 
    + R_h w(t_{n+1/2}) - ( R_h w^{n+1}  + R_h w^n)/2 
    \quad\mbox{(note that $R_hw=\partial_tR_h\eta$)} \label{defect-dw}\\
    d_\eta(\xi_h) &= ((R_h w^{n+1} - R_h w^n)/\tau, \xi_h)_{\Omega_{2}}-((\partial_{tt} \eta^n + \partial_{tt} \eta^{n+1})/2, \xi_h)_{\Omega_{2}}\label{defect-dxi}\\
    d_u(v_h) &= ((R_h u^{n+1} - R_h u^n)/\tau, v_h)_{\Omega_{1}} - ((\partial_t u^n + \partial_t u^{n+1})/2, v_h)_{\Omega_{1}}\notag\\
    & \quad\, + ((u^{n+1} + u^n )/2,v_h)_{\Omega_1} - ((R_h u ^{n+1} + R_h u^n )/2, v_h)_{\Omega_1}\notag\\
    & \quad\, + ((f(t_n)+f(t_{n+1})/2, v_h))_{\Omega_1} - (f(t_{n+1/2}),v_h)_{\Omega_1} \notag\\
    & \quad\, + ((g(t_n)+g(t_{n+1})/2, v_h))_{\Omega_1\backslash\Gamma} - (g(t_{n+1/2}),v_h)_{\Omega_1\backslash\Gamma}\label{defect-du}.
\end{align}
\end{subequations}

By subtracting \eqref{defects-eqn} from \eqref{numerical_weak} and defining the error functions \(e_\eta^n = \eta_h^n - R_h\eta^n\), \(e_w^n = w_h^n - R_h w^n\), \(e_u^n = u_h^n - R_hu^n\), and \(e_p^{n+1/2} = p_h^{n+1/2} - R_h p ^{n+1/2}\), we obtain the following error equations:
\begin{subequations}\label{error-eqn}
    \begin{align}
        &(e_\eta^{n+1}-e_\eta^{n})/\tau = (e_w^{n+1}+e_w^{n})/2 - d_w\quad\mbox{on}\,\,\, \Omega_2, \label{error-eqn-1}\\[5pt] 
    &((e_w^{n+1}-e_w^{n})/\tau , \xi_h )_{\Omega_{2}} 
        + (\nabla (e_\eta^{n+1}+e_\eta^{n})/2, \nabla \xi_h)_{\Omega_{2}} 
        + ((e_u^{n+1} - e_u^n)/\tau, v_h)_{\Omega_{1}} \notag\\
        &\quad\, + (D((e_u^{n+1}+e_u^{n})/2 ) , D(v_h))_{\Omega_{1}}  - (e_p^{n+1/2} , \nabla \cdot v_h )_{\Omega_{1}} + (\nabla \cdot (e_u^{n+1}+e_u^{n})/2, q_h)_{\Omega_{1}} \notag\\
        & = -d_\eta(\xi_h) - d_u(v_h) . \label{error-eqn-2}  
    \end{align}
\end{subequations}
Note that $e_w^{n+1} = e_u^{n+1}$ on $\Gamma$, and \eqref{error-eqn-2}  holds for all test functions $(\xi_h, v_h, q_h) \in X_h(\Omega_{2}) \times X_h(\Omega_{1})\times M_h(\Omega_{1})$ such that $\xi_h = v_h$ on interface $\Gamma$.

\subsection{Defect estimates}
In this section, we will prove the following estimates for the defect terms. 
\begin{lemma}\label{Lemma:defect-estimates}
Under Assumptions \ref{AS0}--\ref{AS2}, the following estimates hold: 
    \begin{subequations}\label{defect-estimates}
        \begin{align}
             \|\nabla d_w\|_{L^\infty(0,T;L^2(\Omega_2))}  &\le C \tau^2,\label{d_w}\\
        |d_\eta(\xi_h)| &\le C(\tau^2+h^{k+1}) \|\xi_h\|_{L^2(\Omega_2)},\label{d_eta}\\
        |d_u(v_h)| &\le C(\tau^2+h^{k+1}) \|v_h\|_{L^2(\Omega_1)}.\label{d_u}
        \end{align}
    \end{subequations} 
\end{lemma}

\begin{proof}
Since the exact solution is sufficiently smooth, from \eqref{diff_H1_u_triple} we can obtain the following result by applying the triangle inequality: 
              \begin{align}\label{ttt-nabla-eta-u}
                \|\partial_{ttt}\nabla R_h \eta(t)\|_{L^\infty(0,T; L^2(\Omega_2))}   \le C.
              \end{align}
In the expression of $d_w$ in \eqref{defect-dw}, by using Taylor's formula and \eqref{ttt-nabla-eta-u}, as well as relation $R_h w = \partial_t R_h \eta$, we have
              \begin{align*}
                \|\nabla d_w\|_{L^\infty(0,T;L^2(\Omega_2))} &\le  \|(\nabla R_h \eta^{n+1} - \nabla R_h\eta^n)/\tau - \partial_t \nabla R_h \eta(t_{n+1/2})\|_{L^\infty(0,T;L^2(\Omega_2))} \\
                &\quad\,+ \|\nabla R_h w(t_{n+1/2}) - (\nabla R_h w^{n+1}  +\nabla R_h w^n)/2\|_{L^\infty(0,T;L^2(\Omega_2))}\\
                & \le C\tau^2 . 
              \end{align*}
This proves \eqref{d_w}. 
%

Substituting relation $R_hw=\partial_tR_h\eta$ into the expression of $d_\eta(\xi_h)$ in \eqref{defect-dxi} and using the triangle inequality, we have 
              \begin{align}\label{rewrite-eta}
                |d_\eta(\xi_h)| & \le |((\partial_t R_h \eta^{n+1} - \partial_t R_h \eta^n)/\tau, \xi_h)_{\Omega_{2}}  -((\partial_t \eta^{n+1} - \partial_t \eta^n)/\tau, \xi_h)_{\Omega_{2}} | \notag\\
                &\quad \,+ |((\partial_t \eta^{n+1} - \partial_t \eta^n)/\tau, \xi_h)_{\Omega_{2}} - ((\partial_{tt} \eta^{n+1} +\partial_{tt} \eta^n)/2, \xi_h) _{\Omega_{2}}| \notag\\
                & =:d_{\eta,1}(\xi_h) + d_{\eta,2}(\xi_h) ,
              \end{align}
with  
\begin{align*}
d_{\eta,1}(\xi_h) 
& =  \Big|\Big( \frac{(\partial_t R_h \eta^{n+1} -\partial_t \eta^{n+1} ) - (\partial_t R_h \eta^n -\partial_t \eta^n)}{\tau} , \xi_h \Big)_{\Omega_{2}} \Big | \notag\\
&\le \sup_{t\in[0,T]} \| \partial_t(\partial_t R_h \eta -\partial_t \eta) \|_{L^2(\Omega_2)} \|\xi_h\|_{L^2(\Omega_2)} \\
&= \sup_{t\in[0,T]} \| \partial_{tt}( R_h \eta - \eta) \|_{L^2(\Omega_2)} \|\xi_h\|_{L^2(\Omega_2)} 
\le Ch^{k+1} \|\xi_h\|_{L^2(\Omega_2)}  ,
\end{align*} 
where the last inequality follows from the second result of Theorem \ref{THM-Ritz}. Moreover, by applying Taylor's formula directly to the expression of $d_{\eta,2}(\xi_h)$, one can derive the following result:  
$$d_{\eta,2}(\xi_h)\le C\tau^2 \sup_{t\in[0,T]} \| \partial_{tttt}\eta\|_{L^2(\Omega_2)} \|\xi_h\|_{L^2(\Omega_2)} 
\le C\tau^2 \sup_{t\in[0,T]} \|\xi_h\|_{L^2(\Omega_2)}  . $$ 
This proves \eqref{d_eta}. The proof of \eqref{d_u} is similar and therefore omitted. 
\end{proof}

\subsection{Error estimates}\label{section:error-estimates}

Finally, by choosing $(\xi_h, v_h, q_h)=((e_w^{n+1}   +e_w^n)/2,(e_u^{n+1}   +e_u^n)/2,e_p^{n+1/2}) \in X_h(\Omega_{2}) \times X_h(\Omega_{1})\times M_h(\Omega_{1})$ in error equation \eqref{error-eqn}, we obtain the following result: 
\begin{align}
    & \frac{\|e_w^{n+1}\|_{L^2(\Omega_2)}^2 -\|e_w^{n}\|_{L^2(\Omega_2)}^2}{2\tau} + \frac{\|\nabla e_\eta^{n+1}\|_{L^2(\Omega_2)}^2 -\|\nabla e_\eta^{n}\|_{L^2(\Omega_2)}^2}{2\tau} \notag\\
    &\quad\, + \frac{\|e_u^{n+1}\|_{L^2(\Omega_1)}^2 -\|e_u^{n}\|_{L^2(\Omega_1)}^2}{2\tau} 
    + \|D((e_u^{n+1} + e_u^n)/2)\|_{L^2(\Omega_1)}^2 \notag\\
    &\le \big|(\nabla(e_\eta^{n+1}+e_\eta^{n})/2, \nabla d_w)_{\Omega_2}\big| + \big|d_\eta((e_w^{n+1}   +e_w^n)/2)\big| + \big|d_u((e_u^{n+1}   +e_u^n)/2)\big|.\notag
\end{align}
Then by using estimates \eqref{defect-estimates} in Lemma \ref{Lemma:defect-estimates}, we obtain
\begin{align}\label{pre-final}
    &\frac{\|e_w^{n+1}\|_{L^2(\Omega_2)}^2 -\|e_w^{n}\|_{L^2(\Omega_2)}^2}{2\tau} + \frac{\|\nabla e_\eta^{n+1}\|_{L^2(\Omega_2)}^2 -\|\nabla e_\eta^{n}\|_{L^2(\Omega_2)}^2}{2\tau}\notag\\
    & \quad\, + \frac{\|e_u^{n+1}\|_{L^2(\Omega_1)}^2 -\|e_u^{n}\|_{L^2(\Omega_1)}^2}{2\tau} 
    + \|D((e_u^{n+1} + e_u^n)/2)\|_{L^2(\Omega_1)}^2\notag \\
    &\le C\epsilon^{-1}(\tau^4 + h^{2k+2} )+ \epsilon( \|\nabla e_\eta^{n+1}\|_{L^2(\Omega_2)}^2 + \|\nabla e_\eta^{n}\|_{L^2(\Omega_2)}^2) \notag\\
    &\quad\,
    + \epsilon (\|e_w^{n+1}\|_{L^2(\Omega_2)}^2 + \|e_w^{n}\|_{L^2(\Omega_2)}^2 )+\epsilon (\|e_u^{n+1}\|_{L^2(\Omega_1)}^2 + \|e_u^{n}\|_{L^2(\Omega_1)}^2) ,
\end{align} 
where $\epsilon$ can be arbitrarily small at the expense of enlarging the term $C\epsilon^{-1}(\tau^4 + h^{2k+2} )$. 

Since \eqref{def-Rh0} can be obtained by choosing $v_h=0$ and $q_h=0$ in \eqref{weak_Ritz} with $\xi_h=0$ on $\Gamma$, it follows that $R_h^0\eta^0 = (R_h\eta)(0)$. Therefore, the initial conditions in \eqref{numerical_initial} guarantees that $e_\eta^0 = 0$, $e_w^0= \partial_t(I_h\eta(0)-R_h\eta(0))$ and $e_u^0=I_hu(0)-R_hu(0)$, and the estimates in Theorem \ref{THM-Ritz} imply that 
$$
\| e_w^0\|_{L^2(\Omega_2)} +\| \nabla e_\eta^0\|_{L^2(\Omega_2)} +  \| e_u^0\|_{L^2(\Omega_1)} \le C(\tau^2 +  h^{k+1}) .
$$
Then, by applying Gr\"onwall's inequality to \eqref{pre-final}, we obtain 
\begin{align}
   \sup_{1\le n\le N} (\| e_w^n\|_{L^2(\Omega_2)} +\| \nabla e_\eta^n\|_{L^2(\Omega_2)} +  \| e_u^n\|_{L^2(\Omega_1)} ) \le C(\tau^2 +  h^{k+1}).
\end{align}
This completes the proof of Theorem \ref{THM-Error}.
\hfill\qed

\section{Numerical examples}

In this section, we present numerical results to support our theoretical analysis on the convergence rates of finite element solutions to the FSI problem. 
 All the computations are performed by the finite element software package NGSolve, which is available at \href{https://ngsolve.org/}{https://ngsolve.org/}.

\begin{example}[Fluid flow in a channel, under periodic boundary condition]\label{section:bloodflow} \upshape

We consider fluid flow in an elastic channel, described by the Stokes equations in the fluid region $\Omega_f = [0,4] \times (1/4, 3/4)$ and the wave equation in two solid regions $\Omega_{s1}=[0,4] \times (3/4, 1)$ and $\Omega_{s2}=[0,4] \times (0, 1/4)$, coupled on the interfaces $\Gamma_1 = [0,4] \times \{3/4\}$ and $\Gamma_2 = [0,4] \times \{1/4\}$. 
%
%
%
%
%
The following exact solution of the FSI problem is constructed (with $\gamma=0.01$):
\[
\begin{aligned}
u(t,x) &= \gamma e^t (2\pi \cos(2\pi x) [-2y^2 + 2y - 3/8], \sin(2\pi x) [4y - 2]) && \text{in } \Omega_f, \\
p(t,x) &= 4\gamma e^t \sin(2\pi x) && \text{in } \Omega_f, \\
\eta_1(t,x) &= \gamma e^t (0, \sin(2\pi x)) && \text{in } \Omega_{s1}, \\
\eta_2(t,x) &= \gamma e^t (0, -\sin(2\pi x)) && \text{in } \Omega_{s2}, 
\end{aligned}
\]
which satisfies the homogeneous Neumann boundary condition on the upper and lower boundaries of $\Omega=\Omega_f\cup\Omega_{s1}\cup \Omega_{s2}\cup\Gamma_1\cup\Gamma_2$, and the periodic boundary condition on the left and right boundaries of $\Omega$. In this case, Assumptions \ref{AS0}--\ref{AS2} are satisfied (as discussed in Scenario 1 below Assumption \ref{AS2}). 

We solve the FSI problem up to $T=0.25$ by the numerical scheme in \eqref{numerical_weak} with P$^{1\rm b}$-P$^1$-P$^{1\rm b}$ and P$^2$-P$^1$-P$^2$ elements for $(u,p,\eta)$, where P$^{1\rm b}$-P$^1$ denotes the Stokes MINI element. The mesh and fluid velocity are illustrated in Figure \ref{mesh_stokes_wave}. 

\begin{figure}[htbp]
    \vspace{-5pt}
    \includegraphics[width=1.0\textwidth]{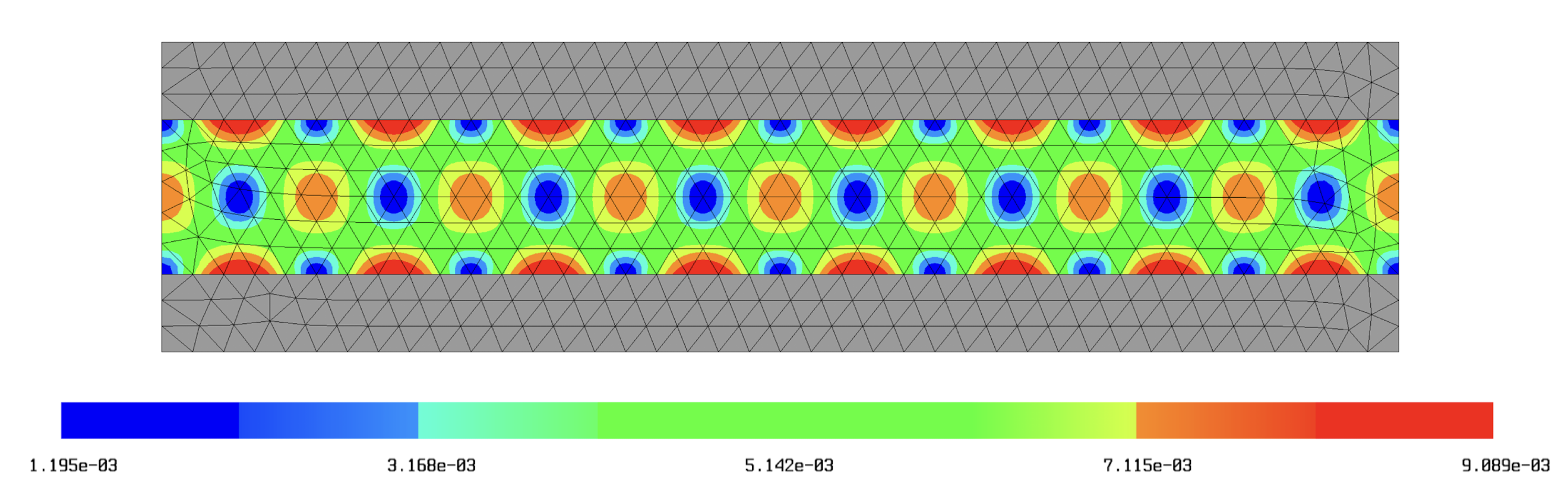}
    \vspace{-10pt}
    \caption{Fluid velocity magnitude (Example \ref{section:bloodflow}).}
    \label{mesh_stokes_wave}
\end{figure}

The errors from spatial discretizations are tested with mesh sizes \(h = 0.03, 0.02, 0.01, 0.008\), using a sufficiently small time stepsize \(\tau = 10^{-3}\) to ensure that the errors from time discretization are negligible in observing the convergence rates with respect to spatial discretizations. Similarly, the errors from time discretizations are tested with time stepsizes \(\tau = 0.02, 0.015, 0.01, 0.008\), using a sufficiently fine mesh size \(h = 0.005\) to ensure that the errors from spatial discretization are negligible in observing the convergence rates with respect to time discretizations. The numerical results in Figures \ref{fig:err_rate_stokes_wave}--\ref{fig:err_rate_stokes_wave_tau} show that the errors are $O(h^{k+1})$ in space and $O(\tau^2)$ in time. This is consistent with the theoretical result proved in Theorem \ref{THM-Error}.  
 


\begin{figure}[htbp]
    \subfigure[P$^{1\rm b}$-P$^1$-P$^{1\rm b}$ element]{\includegraphics[width=0.45\textwidth]{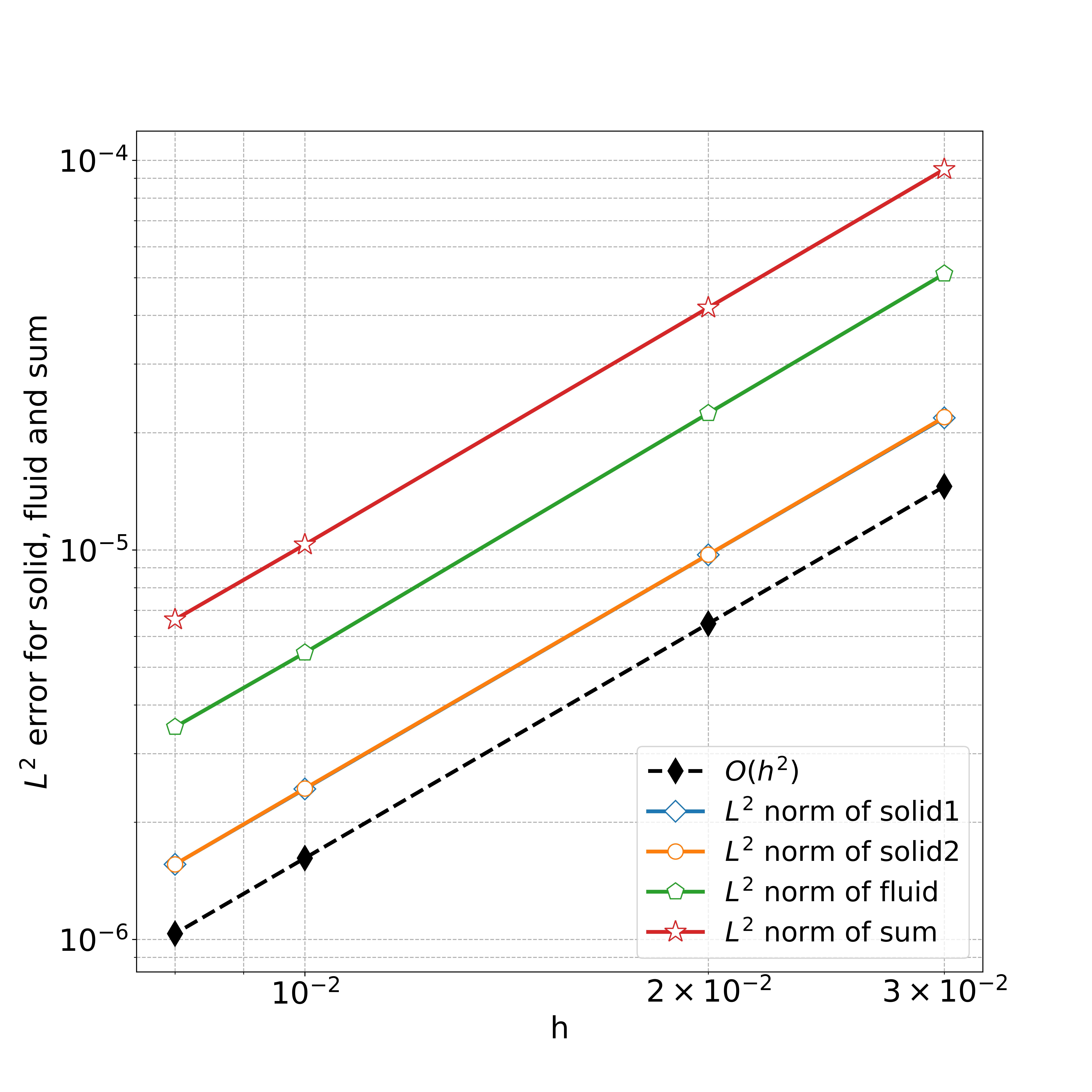}}
    \subfigure[P$^2$-P$^1$-P$^2$ element]{\includegraphics[width=0.45\textwidth]{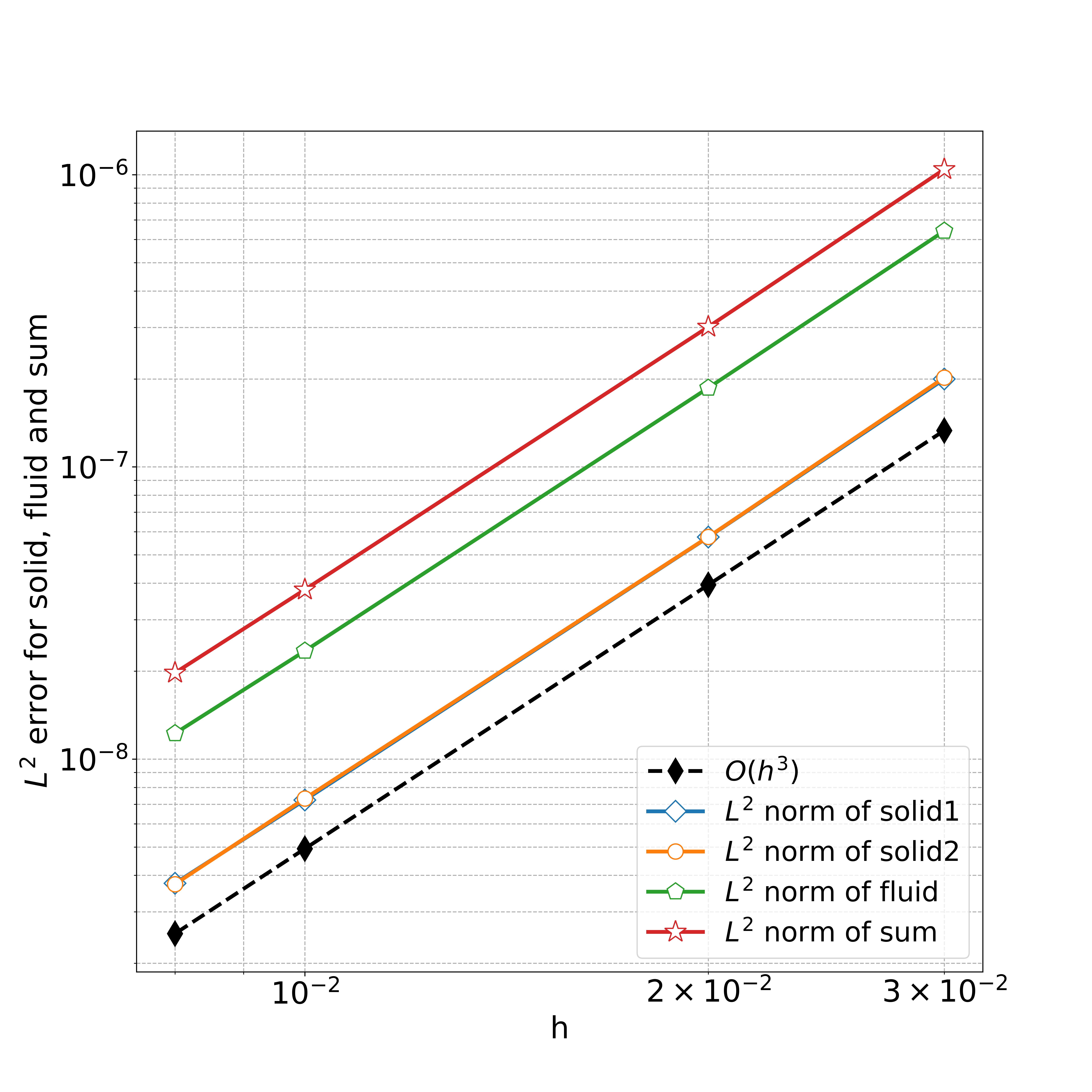}}
    \vspace{-10pt}
    \caption{Spatial discretization errors (Example \ref{section:bloodflow}).}
    \label{fig:err_rate_stokes_wave}
\end{figure}


\begin{figure}[htbp]
    \vspace{-15pt}
    \includegraphics[width=0.45\textwidth]{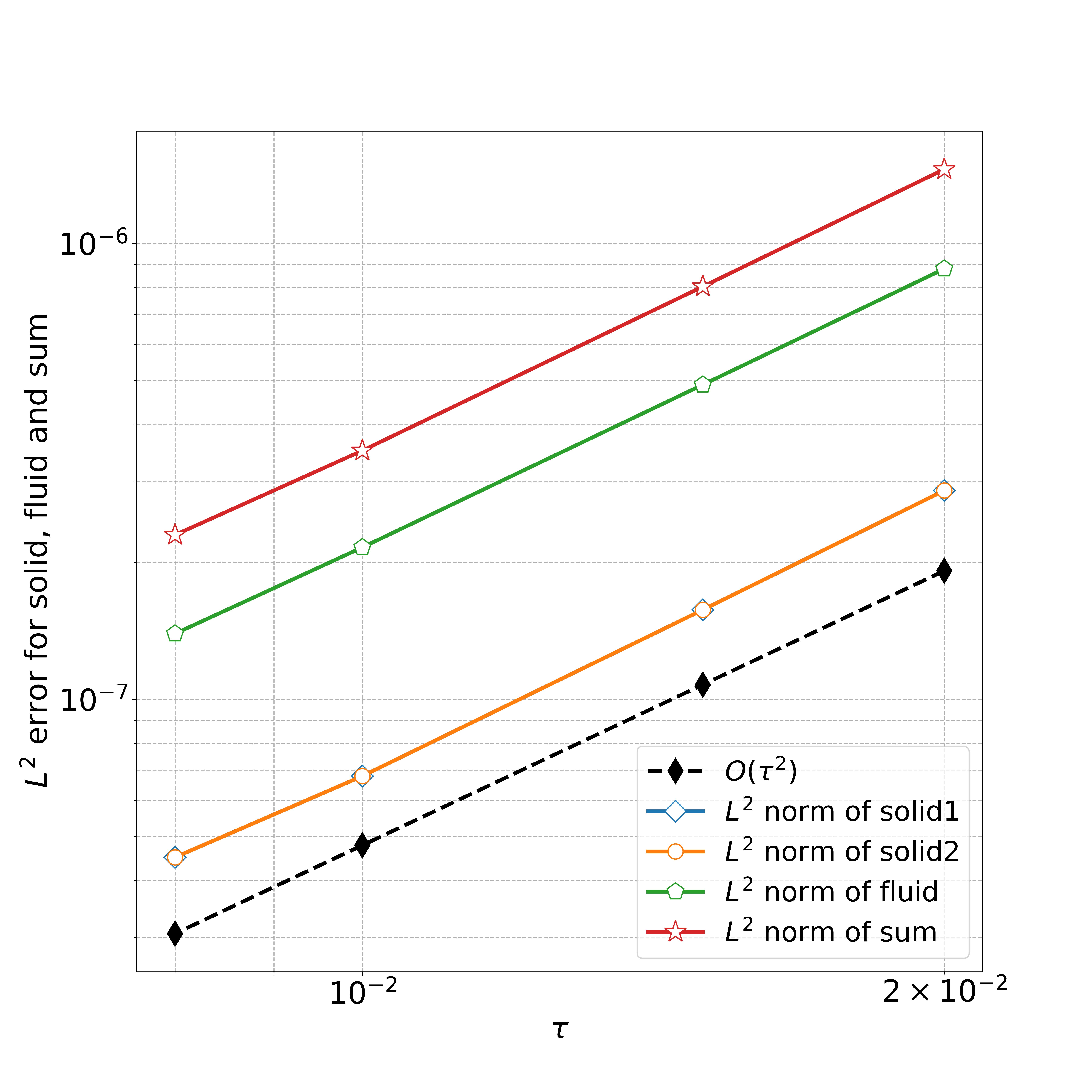}
    \vspace{-10pt}
    \caption{Time discretization errors (Example \ref{section:bloodflow}).}
    \label{fig:err_rate_stokes_wave_tau}
\end{figure}

\end{example}

\begin{example}[Fluid flow in a channel, under the traction boundary condition]\label{section:bloodflow2} \upshape

We consider a similar FSI problem in the domain $\Omega = [0,1]^2=\Omega_f\cup\Omega_{s1}\cup \Omega_{s2}\cup\Gamma_1\cup\Gamma_2$, with \(\Omega_f = [0,4] \times [0.15, 0.85]\), \(\Omega_{s1} = [0,4] \times [0.85, 1]\), \(\Omega_{s2} = [0,4] \times [0, 0.15]\), \(\Gamma_1 = [0,1] \times \{0.85\}\), and \(\Gamma_2 = [0,1] \times \{0.15\}\). The homogeneous Neumann boundary conditions are imposed on the upper and lower boundaries of $\Omega$, the zero traction condition is imposed on the right (outflow) boundary of $\Omega$, and the following traction condition is imposed on the left (inflow) boundary of $\Omega$: 
\[
(Du  - pI) n = - \frac14 [1-\cos(\pi t/2)]^210^4 (y - 0.15)(0.85 - y)n .
\]
In this case, Assumptions \ref{AS0}--\ref{AS2} are satisfied up to an $\epsilon$ modification (the solution is not sufficiently smooth but in $H^{2-\epsilon}$, as discussed in Scenario 2 below Assumption \ref{AS2}). In this case, the convergence rates of the numerical solutions should be $O(\tau^2+h^{2-\epsilon})$, where $\epsilon$ can be arbitrarily small. 

Since the exact solution of this example is unknown, we test the errors of the numerical solutions by using reference solutions computed with sufficiently small mesh size and time stepsize. We solve the FSI problem up to $T=0.1$ by the numerical scheme in \eqref{numerical_weak} with P$^{1\rm b}$-P$^1$-P$^{1\rm b}$ finite element for $(u,p,\eta)$. The numerical solution of fluid velocity magnitude is illustrated in Figure \ref{fig:velocity_stokes_wave}. 
%

\begin{figure}[htbp]
    \includegraphics[width=1.0\textwidth]{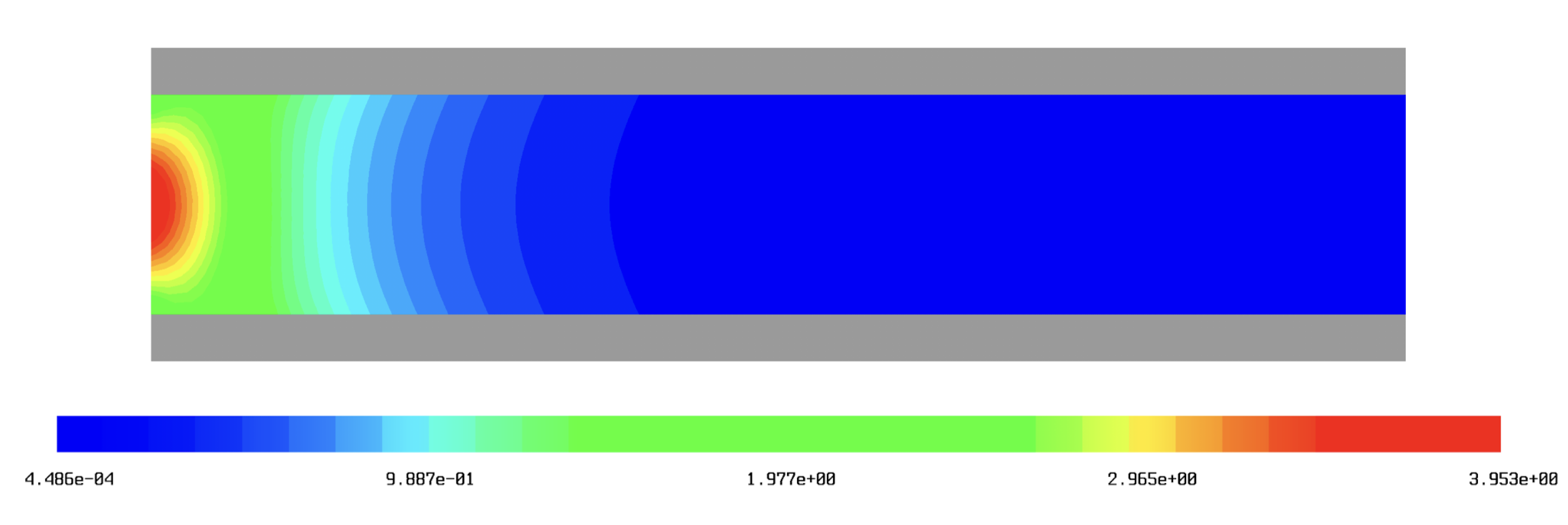}
    \caption{Fluid velocity magnitude (Example \ref{section:bloodflow2}).}
    \label{fig:velocity_stokes_wave}
\end{figure}

The errors from spatial discretizations are tested with a fixed time stepsize $\tau = 0.01$ and mesh sizes $h = 0.03, 0.02, 0.01, 0.008$, using a sufficiently small mesh size \(h = 0.002\) for the reference solution. Similarly, the errors from time discretizations are tested with a fixed mesh size $h = 0.01$ and time stepsizes \(\tau = 0.05, 0.02, 0.01, 0.005\), using a sufficiently small time stepsize $\tau = 0.001$ for the reference solution. The numerical results in Figure \ref{fig:err_rate_stokes_Neumann} show that the errors are about second order in both space and time. This is consistent with the theoretical analysis in this paper.  

%

\begin{figure}[htbp]
    \vspace{-15pt}
    \subfigure[Space discretization errors]{\includegraphics[width=0.45\textwidth]{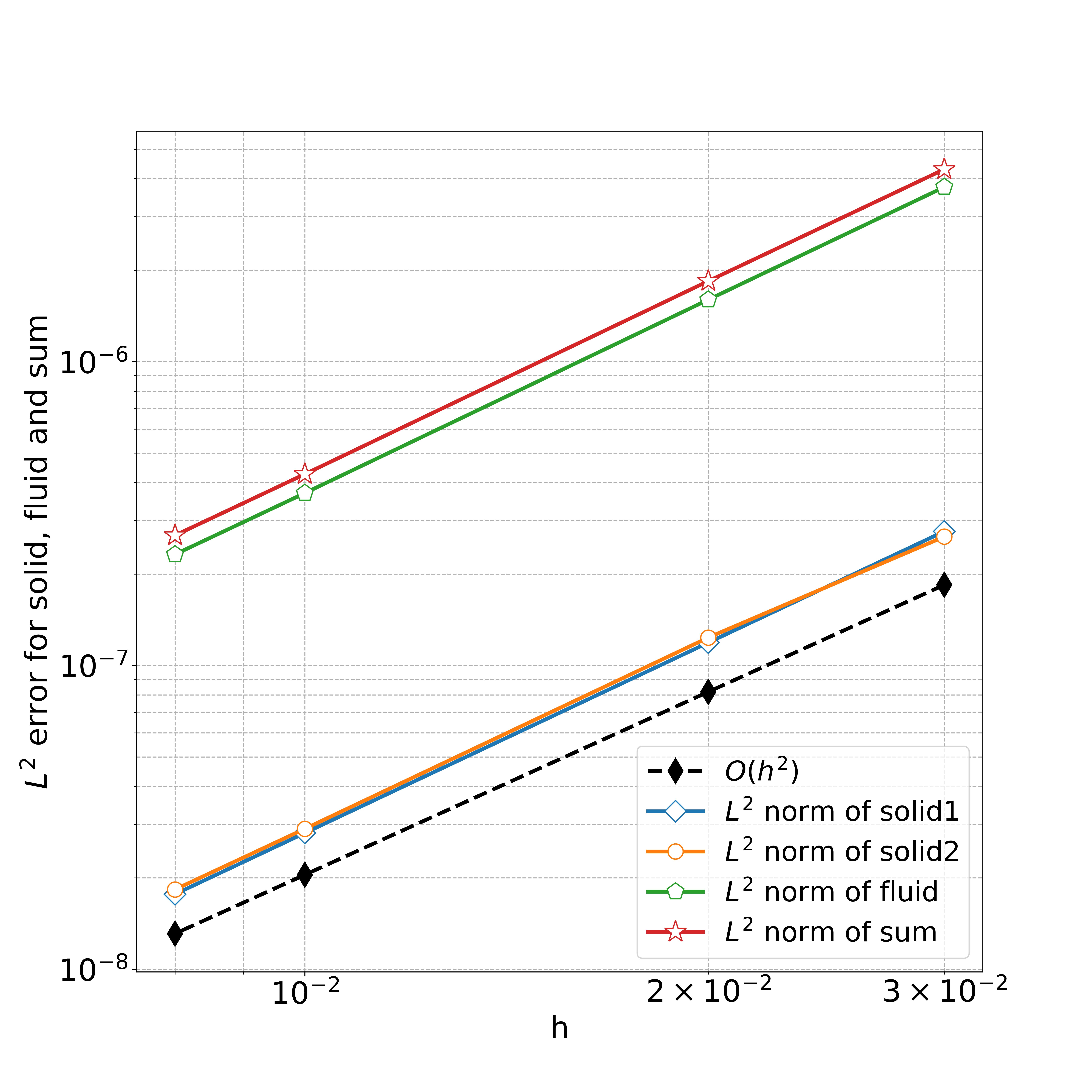}}
    \subfigure[Time discretization errors]{\includegraphics[width=0.45\textwidth]{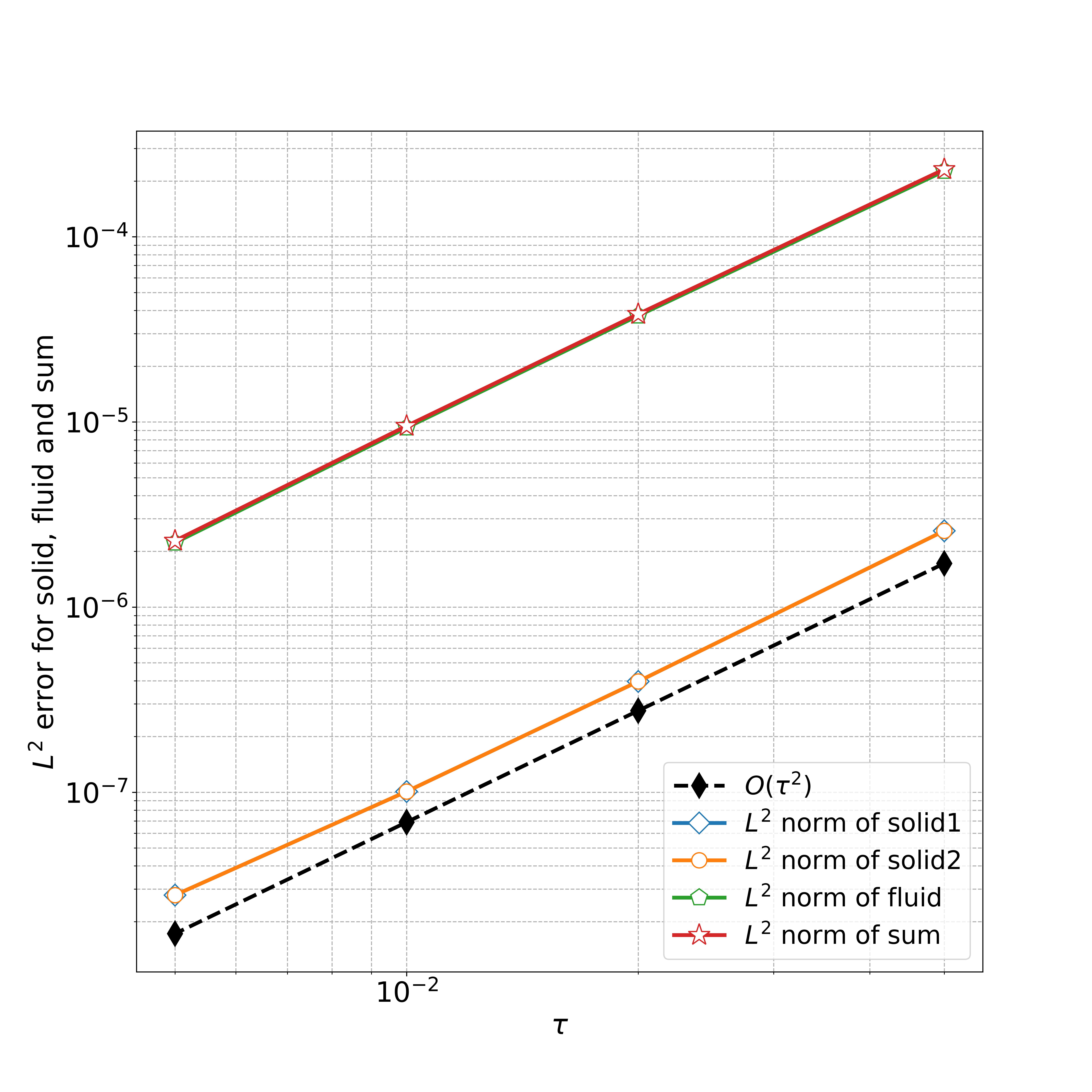}}
    \vspace{-5pt}
    \caption{Spatial errors and convergence rates (Example \ref{section:bloodflow2}).}
    \label{fig:err_rate_stokes_Neumann}
\end{figure}

\end{example}

\begin{example}[The heat-wave system]\label{section:heat} \upshape

We consider an example of heat-wave interaction from \cite{Burman-Fernandez-2023}, described by the heat equation in $\Omega_1=(0,1) \times (0, 3/4)$ and the wave equation in $\Omega_2 = (0,1) \times (3/4,1)$, i.e., 
\begin{equation}\label{para-hyper-para}
    \left\{
    \begin{aligned}
        \partial_t u - \Delta u &= f && \text{in}\,\,\, (0, T) \times \Omega_1 \\
        \partial_{tt}\eta - \Delta \eta &= f && \text{in}\,\,\, (0, T) \times \Omega_2 ,
    \end{aligned}
    \right.
\end{equation}
%
coupled on the interface $\Gamma = [0,1] \times \{3/4\}$: 
\begin{equation}\label{para-hyper-interface}
    \left\{
    \begin{aligned}
        w &= u && \text{on } \Gamma \\
        \nu_s \partial_n \eta &= \nu_f \partial_n u && \text{on } \Gamma . 
    \end{aligned}
    \right.
\end{equation}
The function $f=e^t \left[ \sin(2\pi x) y(1-y) + 2\sin(2\pi x) + 4\pi^2 \sin(2\pi x) y(1-y) \right]$ and the initial conditions are determined by the following exact solution: 
\begin{equation}
    \begin{aligned}
        \eta(t,x) &= e^t \sin(2\pi x) y(1-y), \\
        w(t,x) &= e^t \sin(2\pi x) y(1-y), \\
        u(t,x) &= e^t \sin(2\pi x) y(1-y) . 
    \end{aligned}
\end{equation}
which satisfies the homogeneous Dirichlet boundary conditions on the boundary of $\Omega=(0,1)\times(0,1)$. 
%
%
Our analysis of the dynamic Ritz projection and the optimal-order convergence in $L^\infty(0,T;L^2)$ norm can be extended to the heat-wave interaction problem as well. 
%

We solve this problem up to time $T=0.25$ by a Crank--Nicolson FEM which is similar to \eqref{numerical_weak}, and examine the errors of the numerical solutions in approximating the exact solution. The errors from spatial discretizations are tested with mesh sizes $h = 0.05, 0.04, 0.03, 0.02$, using a sufficiently small time stepsize $\tau = 0.0001$ to ensure that the errors from time discretization is negligible. Similarly, the errors from time discretizations are tested with time stepsizes $\tau = 0.06, 0.05, 0.04, 0.03$, using a sufficiently small mesh size $h = 0.005$. The numerical results are presented in Figure \ref{fig:err_rate_heat_wave_tau}, which shows that the spatial discretization errors and time discretization errors are $O(h^{k+1})$ and $O(\tau^2)$, respectively. This is consistent with the theoretical analysis in this paper.

We test the numerical solution errors for mesh sizes $h = 0.05, 0.04, 0.03, 0.02$ with a sufficiently small time step $\tau = 10^{-4}$ to ensure that spatial discretization errors dominate. The errors, measured in the discrete $L^2$ norms of $e_\eta$ and $e_u$ at $T=0.25$, are shown in Figure \ref{fig:err_rate_heat_wave}. The results confirm that the spatial errors converge at the expected rate $O(h^{k+1})$ for finite elements of degree $k = 1, 2$, consistent with the theoretical error estimates in Theorem \ref{THM-Error}.

\begin{figure}[htbp]
    \vspace{-15pt}
    \subfigure[order $k=1$]{\includegraphics[width=0.45\textwidth]{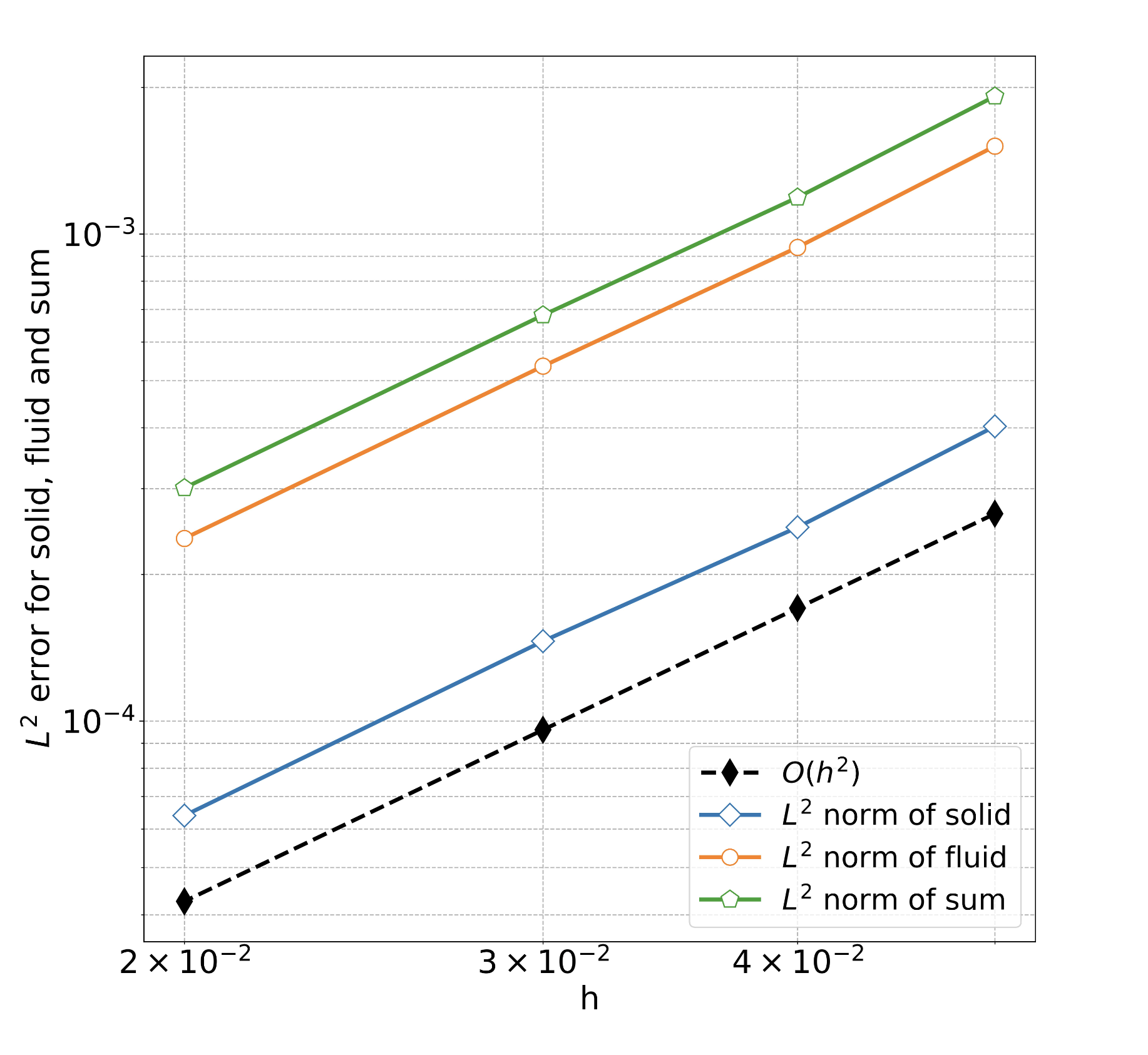}}
    \subfigure[order $k=2$]{\includegraphics[width=0.435\textwidth]{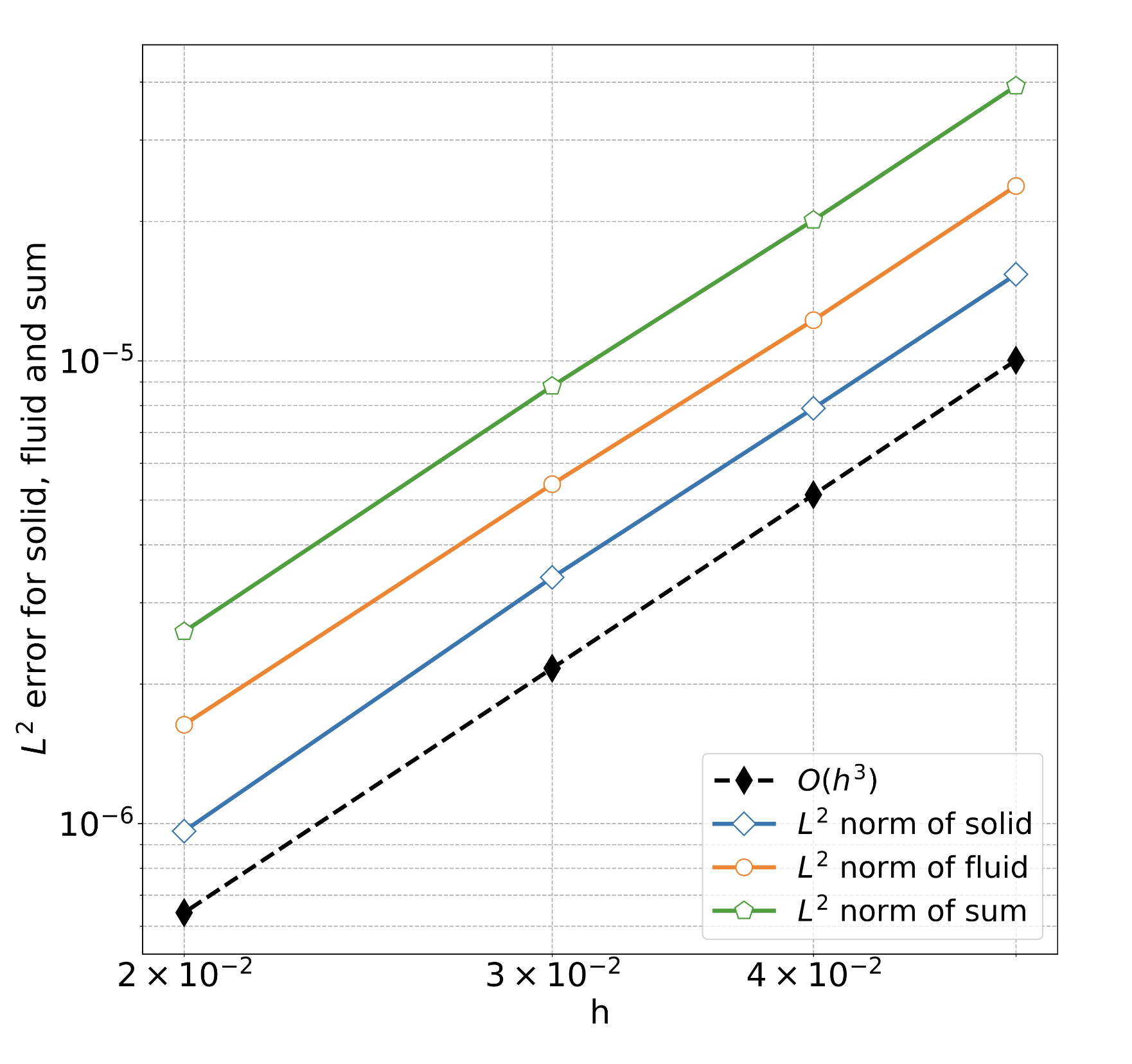}}
    \vspace{-5pt}
    \caption{Spatial discretization errors (Example \ref{section:heat})}
    \label{fig:err_rate_heat_wave}
\end{figure}

\begin{figure}[htbp]
 \vspace{5pt}
    \includegraphics[width=0.45\textwidth]{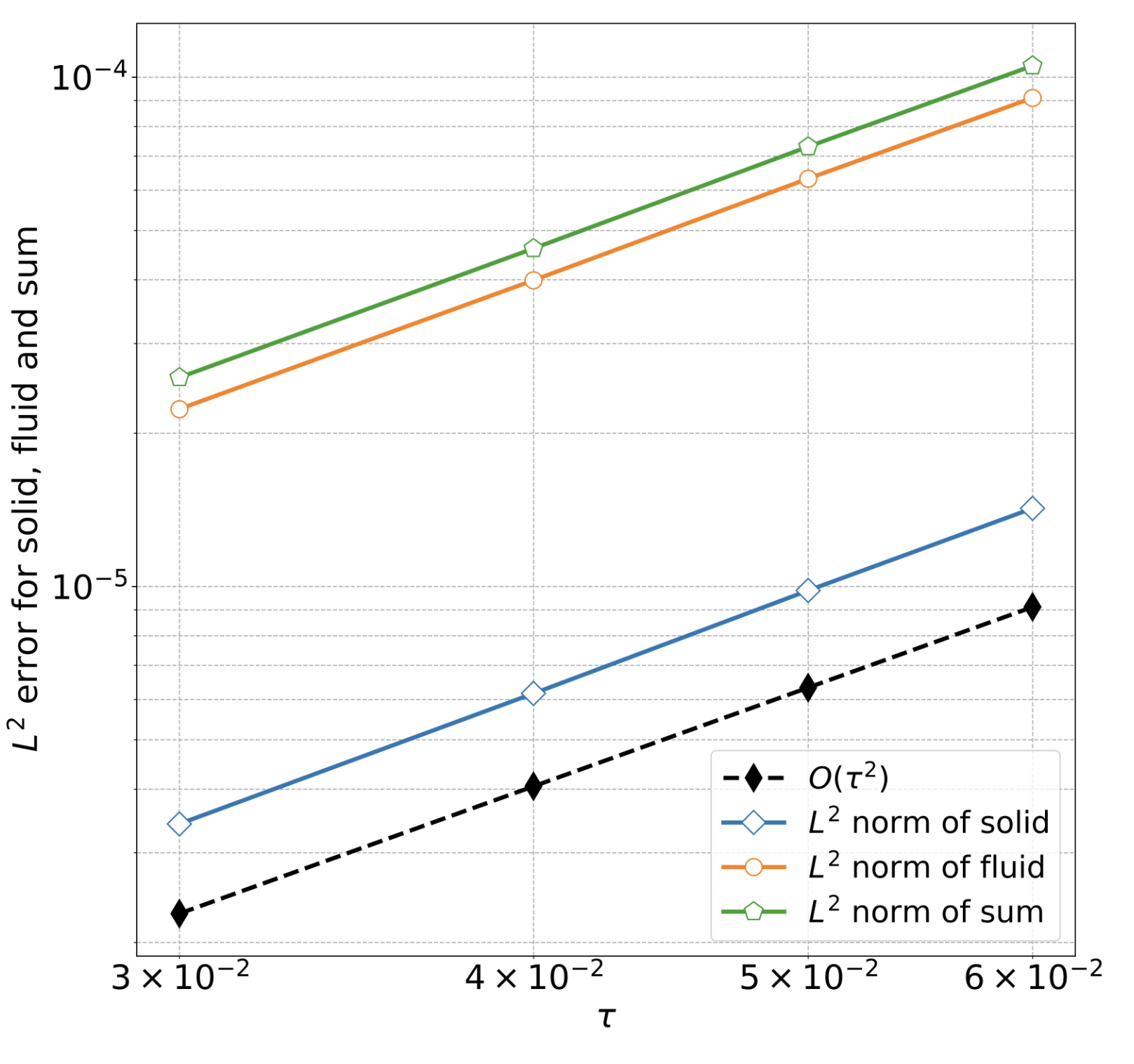}
    \caption{Time discretization errors (Example \ref{section:heat})}
    \label{fig:err_rate_heat_wave_tau}
\end{figure}

\end{example}

\section{Conclusion}

We have defined a finite element dynamic Ritz projection (which satisfies a dynamic interface condition) for the fluid-structure interaction problem in \eqref{exact_eqn}--\eqref{exact_interface}. We have proved existence and uniqueness of the dynamic Ritz projection of the solution, as well as estimates of the error between the solution and its dynamic Ritz projection. By utilizing these results, we have proved optimal-order convergence of finite element solutions to the FSI problem in the $L^\infty(0,T;L^2)$ norm. The result is supported by the numerical examples. 
The analysis of the dynamic Ritz projection and convergence rates of finite element solutions in this paper are based on some $H^2$ regularity assumptions for the Stokes and Poisson equations in the fluid-structure interaction problem; see Assumptions \ref{AS-Ext}--\ref{AS2}. This excludes some cases (such as some boundary conditions and some situations where the interface intersects the boundary) which would lead to corner singularities in the solution. Extension of the error analysis to such cases will be considered in the future.

\end{document}